\documentclass[12pt]{amsart}
\usepackage{amsmath,amssymb,amstext,amsthm,amscd}

\theoremstyle{plain}
\newtheorem{thm}{Theorem}[section]
\newtheorem{pro}[thm]{Proposition}
\newtheorem{lem}[thm]{Lemma}
\newtheorem{cor}[thm]{Corollary}
\theoremstyle{definition}
\newtheorem*{ack}{Acknowledgements}
\theoremstyle{remark}
\newtheorem*{rem}{Remark}

\def\ds{\displaystyle}
\def\pf{{\em Proof.}\ \,}
\def\ov{\overline}
\def\su{\subseteq}
\def\({\left(}
\def\){\right)}

\def\inc{\hookrightarrow}
\def\fun{\rightarrow}
\def\lfun{{\longrightarrow}}

\def\om{\omega}

\def\bC{{\mathbb{C}}}
\def\bP{{\mathbb{P}}}
\def\bZ{{\mathbb{Z}}}

\def\cA{{\mathcal A}}
\def\cB{{\mathcal B}}
\def\BL{{\mathcal B}{\mathcal L}}
\def\cE{{\mathcal E}}
\def\cExt{{\mathcal E}xt}
\def\cF{{\mathcal F}}
\def\cHom{{\mathcal H}om}
\def\cI{{\mathcal I}}
\def\cN{{\mathcal N}}
\def\cO{{\mathcal O}}
\def\SU{{\mathcal S}{\mathcal U}}
\def\cTor{{\mathcal T}or}

\def\wt{\widetilde}
\def\tC{{\widetilde{C}}}
\def\tE{{\widetilde{E}}}
\def\tG{{\widetilde{\Gamma}}}
\def\tH{{\widetilde{H}}}
\def\tX{{\widetilde{X}}}
\def\tY{{\widetilde{Y}}}
\def\tp{{\widetilde{p}}}

\def\codim{{\textrm{codim\,}}}
\def\Ext{{\textrm{Ext\,}}}
\def\eExt{{\emph{Ext}\,}}
\def\Hom{{\textrm{Hom\,}}}
\def\eHom{{\emph{Hom}\,}}
\def\Im{{\textrm{Im\,}}}
\def\eIm{{\emph{Im}\,}}
\def\Pic{{\textrm{Pic}}}
\def\Span{{\textrm{Span\,}}}
\def\Tors{{\textrm{Tors\,}}}

\def\oSU{{\ov{\SU_C(2,L)}}}

\begin{document}

\title{Extensions and rank-$2$ vector bundles
on irreducible nodal curves}

\author{D.\ Arcara}
\address{Department of Mathematics, University of Utah,
155 S. 1400 E., Room 233, Salt Lake City, UT 84112-0090, USA}
\email{arcara@math.utah.edu}

\begin{abstract}
We generalize Bertram's work on rank two vector bundles to an irreducible
projective nodal curve $C$.
We use the natural rational map
$\phi_L\colon\bP(\Ext^1_C(L,\cO_C))\fun\SU_C(2,L)\su\oSU$
defined by $\phi_L([0\fun\cO_C\fun E\fun L\fun0])=E$ to study a
compactification $\oSU$ of the moduli space $\SU_C(2,L)$ of semi-stable vector
bundles of rank $2$ and determinant $L$ on $C$.
In particular, we resolve the indeterminancy of $\phi_L$ in the case
$\deg L=3,4$ via a sequence of three blow-ups with smooth centers.
\end{abstract}

\maketitle

\tableofcontents

\section{Introduction}

In [Ber92], Bertram used extensions of line bundles to study rank-$2$ vector
bundles of fixed determinant on a smooth curve.
We generalize his construction to an irreducible projective nodal curve $C$.
The idea is to consider extensions of $L$ by $\cO_C$, where $L$ is a generic
line bundle on $C$, and consider the `forgetful' map which sends an extension
to the vector bundle of rank $2$ in the middle, forgetting the extension maps.
This gives a rational map from $\bP(\Ext^1_C(L,\cO_C))$ to $\SU_C(2,L)$, the
moduli space of semi-stable vector bundles of rank $2$ and determinant $L$.
If the arithmetic genus of $C$ is  $\geq2$ and $\deg L=3$ or $4$, we resolve
the indeterminancy of the map by a sequence of three blow-ups with smooth
centers.
A nice aspect of these blow-ups is that there exists at each stage a
`universal bundle' which induces the rational map in a natural way.

Let $\oSU$ be the natural compactification of $\SU_C(2,L)$ via torsion-free
sheaves introduced by Newstead and Seshadri (see [New78] and [Ses82]).
Our main theorem is the following.

\begin{thm}
Let $C$ be an irreducible projective nodal curve of arithmetic genus $\ge2$,
and let $L$ be a generic line bundle on $C$ of degree $3$ or $4$.
Let $\phi_L\colon\bP(\eExt^1_C(L,\cO_C))\fun\SU_C(2,L)\su\oSU$ be the natural
rational map defined by $\phi_L([0\fun\cO_C\fun E\fun L\fun0])=E$.
There exist a sequence of three blow-ups with smooth centers
$$ \bP_{L,3} \overset{\varepsilon_3}{\lfun} \bP_{L,2}
\overset{\varepsilon_2}{\lfun} \bP_{L,1} \overset{\varepsilon_1}{\lfun}
\bP(\eExt^1_C(L,\cO_C)), $$
such that
$$ \phi_L \circ \varepsilon_1 \circ \varepsilon_2 \circ \varepsilon_3 \colon
\bP_{L,3} \lfun \oSU $$
extends to a morphism $\phi_{L,3}$.
\end{thm}

An important fact is that the fibers of $\phi_{L,3}$ are connected.
As a corollary, we can give a new proof of the fact that, if $C$ is an
irreducible nodal curve of arithmetic genus $2$, and $\deg L$ is odd, then the
normalization morphism $\oSU^\nu\fun\oSU$ is one-to-one.

We also give the idea for a new proof of $\oSU\simeq\bP^3$ for an
irreducible nodal curve of arithmetic genus $2$ when $\deg L$ is even.

If the arithmetic genus of $C$ is $1$, using the morphism $\phi_L$ with $L$ of
degree $1$ or $2$, we prove that, as in the smooth case,
$$ \oSU = \SU_C(2,L) \simeq
\begin{cases}
\{ pt \} & \textrm{if } \deg L \textrm{ is odd} \\ 
\bP^1 & \textrm{if } \deg L \textrm{ is even}
\end{cases}, $$
and there are no stable bundles of even degree.

In general, using the rational map $\phi_{L,3}$ with $\deg L\geq3$, we prove
that the complement of $\SU_C(2,L)$ in $\oSU$ has codimension $\geq3$ for
every irreducible nodal curve of arithmetic genus $\geq2$.
It follows, using [Bho99] and [Bho04], that $A_{3g-4}(\oSU)\simeq\bZ$.
Moreover, if $\deg L=2g-1$, we find open subsets $U\su\bP_{L,3}$ and
$V\su\SU_C(2,L)$ such that $\phi_{L,3}|_U\colon U\fun V$ is an isomorphism,
and $\codim(\oSU\setminus V,\oSU)\geq2$.
As a corollary, we prove directly that
$$ A_{3g-4}(\SU_C(2,L)) = A_{3g-4}(\oSU) \simeq \bZ $$
if $\deg L$ is odd.

\begin{ack}
A special thanks to my Ph.D.\ adviser E.\ Izadi for her constant guide and
support throughout the program.
I would also like to thank V.\ Alexeev, A.\ Bertram, R.\ Smith,
V.\ Vologodski, and especially R.\ Varley, for many useful discussions.
\end{ack}

\begin{center}
{\bf Notation}
\end{center}

Let $J$ be the set of nodes of $C$.
For a subset $J'$ of $J$, we denote by $\pi_{J'}\colon C_{J'}\fun C$ the
partial normalization of $C$ along the nodes in $J'$.
In particular, for $J'=J$ we obtain the normalization $\pi\colon N\fun C$ of
$C$.
For every $p\in J$, let $p_1,p_2$ be the two points which map to $p$ under any
partial normalization map $\pi_{J'}$ with $p\in J'$.

If $X$ is a projective variety, a sheaf on $X\times C$ of the form
$\pi_X^*F\otimes\pi_C^*G$ (for some sheaves $F$ on $X$ and $G$ on $C$) shall
be denoted by $F\boxtimes G$.
If it is of the form $\pi_X^*F$ or $\pi_C^*G$, we shall sometimes just denote
it by $F$ or $G$, if it is clear from the context that we are actually
considering the sheaf on $X\times C$.

We shall assume throughout the paper, unless it is explicitly stated
otherwise, that the arithmetic genus of $C$ is $g\ge2$ (and we shall simply
call it the genus of $C$).

Whenever we do not explicitly define a homomorphism of extension spaces
throughout this paper, a natural push-forward or pull-back of extensions is
understood.

\section{Description of $\phi_L$}

Let us start with extending to our situation some of the basic results of the
smooth case.
Since $\Ext^1_C(L,\cO_C)\simeq H^1(C,L^{-1})\simeq H^0(C,L\otimes\om_C)^*$
(see [Har77, chapter III]), the linear system $|L\otimes\om_C|$ defines a
rational map $\varphi_{L\otimes\om_C} \colon  C \fun |L\otimes\om_C|^* \simeq
\bP_L$, where we denote $\bP(\Ext^1_C(L,\cO_C))$ by $\bP_L$ to simplify the
notation.
Let $U_L\su\bP_L$ be the open locus of semi-stable extensions, i.e., the open
subset where $\phi_L$ is well-defined.

\begin{pro}\label{UL}
$(1)$ If $\deg L<0$, then $U_L=\emptyset$.

$(2)$ If $0\le\deg L\le2$, then $U_L=\bP_L$.

$(3)$ If $3\le\deg L\le4$, then
$U_L=\bP_L\setminus\varphi_{L\otimes\om_C}(C)$.
\end{pro}

\pf
The same is true for a smooth curve (see [Ber92]), and the proof in our case
is similar, except for two technical details that we prove in Lemmas
\ref{rk1tf} and \ref{emb}.
\qed

\begin{lem}\label{rk1tf}
Every torsion-free sheaf of rank $1$ and degree $1$ on $C$ with a section is
either isomorphic to $\cO_{C}(q)$ for some smooth point $q\in C$ $($if it is
locally-free$)$ or it is isomorphic to $(\pi_p)_*\cO_{C_p}$ for some node $p$
$($if it is not locally-free$)$.
\end{lem}

\pf
Every torsion-free non-locally-free coherent sheaf $F$ of rank $1$ on $C$ is
of the form $(\pi_{J'})_*\cF$ for some line bundle $\cF$ on a partial
normalization $C_{J'}$ of $C$ (see [Ses82]).
If it has a section, then $\cO_C\su F$ implies that
$\cO_{C_{J'}}=\pi_{J'}^*\cO_C\su\cF$, and therefore, $\deg\cF\geq0$.
Since $F$ has degree $1$, and it is not locally-free, $J'$ must contain only
one node $p$, $\cF$ must have degree $0$, and therefore be isomorphic to
$\cO_{C_p}$.
\qed

\begin{lem}\label{emb}
If $q$ is a smooth point of $C$, then
$$ \varphi_{L \otimes \om_C}(q) =
\bP(\ker(\eExt^1_C(L,\cO_C) \overset{\psi_q}{\lfun}\eExt^1_C(L,\cO_C(q)))). $$
If $p$ is a node of $C$, then
$$ \varphi_{L \otimes \om_C}(p) =
\bP(\ker(\eExt^1_C(L,\cO_C) \overset{\psi_p}{\lfun}
\eExt^1_C(L,(\pi_p)_*\cO_{C_p}))). $$
\end{lem}

\pf
If $q$ is a smooth point, the proof is the same as in the smooth case.
If $p$ is a node, then $\varphi_{L\otimes\om_C}(p)$ is the hyperplane of
$H^0(C,L\otimes\om_C)$ defined by the sections vanishing at $p$.
Since the sheaf generated by the regular functions vanishing at $p$ is the
sheaf $(\pi_p)_*(\cO_{C_p}(-p_1-p_2))$ and its dual is $(\pi_p)_*\cO_{C_p}$,
$\varphi_{L\otimes\om_C}(p)$ corresponds to the kernel of the linear
homomorphism $H^1(C,L^{-1})\fun H^1(C,L^{-1}\otimes(\pi_p)_*\cO_{C_p})$, where
we identified $H^0(C,L\otimes\om_C)^*$ with $H^1(C,L^{-1})$.
If $G$ is any coherent sheaf, we can identify $\Ext^1_C(L,G)$ with
$H^1(C,L^{-1}\otimes G)$, and the
linear homomorphism above becomes $\psi_p$ as claimed.
\qed

From now on, all through Section \ref{fibers}, we shall restrict ourselves to
the case when $\deg L$ is either $3$ or $4$.

\begin{lem}
If $\deg L\geq3$, then $\varphi_{L\otimes\om_C}$ is an embedding.
\end{lem}

\begin{rem}
Since $\varphi_{L\otimes\om_C}$ is an isomorphism onto its image, we
shall identify $C$ with $\varphi_{L\otimes\om_C}(C)\su\bP_L$.
\end{rem}

\pf
We need to prove that, for every $q,q'\in C$ not both equal to a node $p$,
$ H^1(C, L \otimes \om_C \otimes \cI_q \otimes \cI_{q'}) =
H^1(C, L \otimes \om_C) = 0, $
where $\cI_q$ [resp.\ $\cI_{q'}$] is the ideal sheaf of the point $q$
[resp.\ $q'$], and that
$ H^1(C, L \otimes \om_C \otimes \cI_p^2) =
H^1(C, L \otimes \om_C) = 0 $
for every node $p$ (see [Bar87]).

Case I: $q,q'$ smooth points.
Then, by Serre duality,
$ H^1(C, L \otimes \om_C \otimes \cI_q \otimes \cI_{q'}) \simeq
H^0(C, L^{-1}(q + q'))^*, $
which is zero because $\deg(L^{-1}(q+q'))<0$.

Case II: $q$ smooth point and $q'=p$ node.
Since $\pi_p^*\om_C\simeq\om_{C_p}(p_1+p_2)$ (see [Bar87]), using the
projection formula we obtain $\om_C\otimes\cI_p\simeq(\pi_p)_*\om_{C_p}$, and
$ L \otimes \om_C \otimes \cI_q \otimes \cI_p \simeq L(-q) \otimes
(\pi_p)_*\om_{C_p} \simeq (\pi_p)_*(\pi_p^*(L(-q)) \otimes \om_{C_p}). $
Therefore,
$ H^1(C, L \otimes \om_C \otimes \cI_q \otimes \cI_{q'}) \simeq
H^1(C_p, \pi_p^*(L(-q)) \otimes \om_{C_p}) \simeq
H^0(C_p, \pi_p^*(L^{-1}(q)))^*, $
which is zero because $\deg(\pi_p^*(L^{-1}(q)))<0$.

Case III: $q=p$, $q'=p'$ distinct nodes.
Since $\pi_{\{p,p'\}}^*\om_C\simeq\om_{C_{\{p,p'\}}}(p_1+p_2+p'_1+p'_2)$ (see
[Bar87]), we obtain
$\om_C\otimes\cI_p\otimes\cI_{p'}\simeq(\pi_{\{p,p'\}})_*\om_{C_{\{p,p'\}}}$,
and
$ L \otimes \om_C \otimes \cI_p \otimes \cI_{p'} \simeq L \otimes
(\pi_{\{p,p'\}})_*\om_{C_{\{p,p'\}}} \simeq 
(\pi_{\{p,p'\}})_*(\pi_{\{p,p'\}}^*L \otimes \om_{C_{\{p,p'\}}}). $
Therefore,
$ H^1(C, L \otimes \om_C \otimes \cI_q \otimes \cI_{q'}) \simeq
H^1(C_{\{p,p'\}}, \pi_{\{p,p'\}}^* L \otimes \om_{C_{\{p,p'\}}}) \simeq
H^0(C_{\{p,p'\}}, \pi_{\{p,p'\}}^* L^{-1})^*, $
which is zero because $\deg(\pi_{\{p,p'\}}^*L^{-1})<0$.

Case IV: $q=q'=p$ node.
Then $\om_C\otimes\cI_p^2\simeq(\pi_p)_*(\om_{C_p}(-p_1-p_2))$ and
$ L \otimes \om_C \otimes \cI_p^2 \simeq
(\pi_p)_*(\pi_p^*L \otimes \om_{C_p}(- p_1 - p_2)). $
Therefore,
$ H^1(C, L \otimes \om_C \otimes \cI_p^2) \simeq
H^1(C_p, \pi_p^*L \otimes \om_{C_p}(- p_1 - p_2)) \simeq
H^0(C_p, \pi_p^*L^{-1} (p_1 + p_2))^*, $
which is zero because $\deg(\pi_p^*L^{-1}(p_1+p_2))<0$.
\qed

\begin{lem}\label{TpC}
The projective tangent plane to $C$ at a node $p$ is
$$ T_p C =
\bP(\ker(\eExt^1_C(L,\cO_C) \overset{\psi_{T_pC}}{\lfun}
\eExt^1_C(L,(\pi_p)_*(\cO_{C_p}(p_1+p_2))))). $$
\end{lem}

\pf
It is easy to see that all the kernels involved in this proof have the right
dimension.
The secant line between the node $p$ and a smooth point $q$ is given by
$$ \bP(\ker(\Ext^1_C(L,\cO_C) \lfun \Ext^1_C(L,(\pi_p)_*\cO_{C_p}(q)))), $$
this being a $1$-dimensional linear subspace of $\bP_L$ which contains
both $p$ and $q$.
If we take the limit as $q\mapsto p$ along the branch corresponding to
$p_i$ ($i=1,2$), we see that the projective tangent line at $p$ to that
branch is
$ X_{p_i} := \bP(\ker(\Ext^1_C(L,\cO_C) \fun
\Ext^1_C(L,(\pi_p)_*(\cO_{C_p}(p_i)))))\ (i = 1,2). $
Since $\bP(\ker(\psi_{T_pC}))$ is a $2$-dimensional linear subspace of
$\bP_L$ which contains both $X_{p_1}$ and $X_{p_2}$, it is the
projective tangent plane $T_pC$ to $C$ at $p$.
\qed

We end this section with an important way to describe the rational map
$\phi_L$.

\begin{pro}\label{cEL}
There exists a locally-free sheaf $\cE_L$ on $\bP_L\times C$ such that
$ \cE_L|_{\{ x \} \times C} \simeq \phi_L(x) $
for every $x\in\bP_L\setminus C$.
Moreover, $\cE_L$ is an extension in
$\eExt^1_{\bP_L\times C}(L,\cO_{\bP_L}(1))$, and for every $a\ne0$ in
$\eExt^1_C(L,\cO_C)$, if we identify
$\pi_{\bP_L}^*\cO_{\bP_L}(1))|_{\{[a]\}\times C}$ with $\cO_C$ using
$a$, $\cE_L$ restricts on $\{[a]\}\times C$ to the extension $a$
itself.
\end{pro}

\pf
Let $\cE_L$ be the extension corresponding to the identity homomorphism under
the natural isomorphism 
$\Ext^1_{\bP_L\times C}(L,\cO_{\bP_L}(1))\simeq
\Hom(\Ext^1_C(L,\cO_C),\Ext^1_C(L,\cO_C))$.
Then, if $a\neq0$ is an extension $0\fun\cO_C\fun E_a\fun L\fun0$,
$\cE_L|_{\{[a]\}\times C}$ is $E_a$ (see [Arc04]).
\qed

\section{The first blow-up}

Since the indeterminancy locus of the rational map
$\phi_L\colon\bP_L\fun\oSU$ is the curve $C\su\bP_L$, to resolve the
indeterminancy via a sequence of blow-ups with smooth centers, we need
to begin the process with the blow-up of $\bP_L$ at the set of nodes
$J\su C$.
By Lemma \ref{emb}, a node $p$ is $\bP(\ker(\psi_p))$, where $\psi_p$ is the
natural linear homomorphism
$\Ext^1_C(L,\cO_C)\fun\Ext^1_C(L,(\pi_p)_*\cO_{C_p})$.
Therefore, the exceptional divisor $E_1$ of
$ \bP_{L,1} := \BL_J \bP_L \overset{\varepsilon_1}{\fun} \bP_L $
is canonically isomorphic to
$\coprod_{p\in J}\bP(\Ext^1_C(L, (\pi_p)_*\cO_{C_p}))$.

\begin{thm}\label{1}
$(a)$ The composition
$ \phi_L \circ \varepsilon_1 \colon \bP_{L,1} \fun \oSU $
extends to a rational map $\phi_{L,1}$ defined as follows:
For every node $p$, a point $x\in E_1|_p$ corresponds to an extension $E'_x$
in $\eExt^1_C(L, (\pi_p)_*\cO_{C_p})$.
Its image $\phi_{L,1}(x)$ is the torsion-free sheaf $E_x$ which is the image
of $E'_x$ under the natural homomorphism
$$ \eExt^1_C(L, (\pi_p)_*\cO_{C_p}) \overset{\psi_{L_p}}{\fun}
\eExt^1_C(L\otimes((\pi_p)_*\cO_{C_p})^*,(\pi_p)_*\cO_{C_p}). $$

$(b)$ The indeterminancy locus of the rational map
$\phi_{L,1}\colon\bP_{L,1}\fun\oSU$ is the union of the strict
transform $\tC_1$ of $C$ and the lines
$ L_p := \bP(\ker(\psi_{L_p})) \su E_1|_p. $
\end{thm}

The strict transform $\tC_1$ of $C$ is isomorphic to $N$ and, for each node
$p$, it intersects $E_1|_p$ at the two points $p_1,p_2$ lying on $p$.
The following lemma describes the lines $L_p$.

\begin{lem}\label{L1}
The points on $L_p$ correspond to the directions tangent to $p$ in $T_pC$, the
projective tangent plane to $C$ at $p$.
In particular, $L_p$ is the line through $p_1$ and $p_2$ in $E_1|_p$.
\end{lem}

\pf
It suffices to show that $L_p$ contains $p_1$ and $p_2$.
It is easy to see that, for $i=1,2$,
$p_i=\bP(\ker(\psi_i))\su\bP(\Ext^1_C(L,(\pi_p)_*\cO_{C_p}))\simeq E_1|_p$,
where $\psi_i$ is the natural map
$\Ext^1_C(L,(\pi_p)_*\cO_{C_p})\fun\Ext^1_C(L,(\pi_p)_*\cO_{C_p}(p_i))$.
To prove that $p_1,p_2\in L_p$, we need to show that
$\ker\psi_i\su\ker\psi_{L_p}$ for $i=1,2$.
A non-trivial extension $E$ in $\Ext^1_C(L,(\pi_p)_*\cO_{C_p})$ is in the
kernel of $\psi_i$ if and only if there exists a surjective map
$E\fun(\pi_p)_*\cO_{C_p}(p_i)$.
The kernel of this map is $L\otimes((\pi_p)_*\cO_{C_p})^*$, and $E$ is
therefore also in the kernel of $\psi_{L_p}$.
\qed

It can be shown that, for every node $p$, the image of the natural linear
homomorphisms $\psi_{L_p}$ is isomorphic to the space of extensions
$\Ext^1_{C_p}(\pi_p^*L(-p_1-p_2),\cO_{C_p})$ via the homomorphism $(\pi_p)_*$.
In particular, no torsion-free sheaf in the image of $\phi_{L,1}|_{E_1}$ is
locally-free, being a push-forward from a partial normalization of $C$.

\begin{cor}
The image $\phi_{L,1}(\bP_{L,1}\setminus(\coprod_{p\in J}L_p\cup\tC_1))$ of
$\phi_{L,1}$ in $\oSU$ is given by\footnote{Note that, since
$((\pi_p)_*\cO_{C_p})^*\simeq(\pi_p)_*\cO_{C_p}(-p_1-p_2)$,
the push-forward of $\pi_p^*L(-p_1-p_2)$ from $C_p$ to $C$ is isomorphic to
$L\otimes((\pi_p)_*\cO_{C_p})^*$ by the projection formula.}
$$ \phi_L(\bP_L \setminus C) \cup \{ E \in \oSU \mid E = (\pi_p)_* \cE
\textrm{ for some } p \in J \textrm{ and } \cE \in
\eExt^1_{C_p}(\pi_p^*L(- p_1 - p_2), \cO_{C_p}) \}. $$
\end{cor}

Before we prove Theorem \ref{1}, we need the following lemma.

\begin{lem}\label{extss}
For every node $p$, all non-trivial extensions
$$ 0 \lfun (\pi_p)_*\cO_{C_p} \lfun E \lfun L \otimes ((\pi_p)_*\cO_{C_p})^*
\lfun 0 $$
in $\eExt^1_C(L\otimes((\pi_p)_*\cO_{C_p})^*,(\pi_p)_*\cO_{C_p})$ are
semi-stable.
\end{lem}

\pf
Assume that $E$ is not semi-stable.
Then there exists a torsion-free quotient $F$ of $E$ of rank $1$ and degree
$\leq1$.
Consider the composite map $(\pi_p)_*\cO_{C_p}\inc E\fun F$.
If it is the zero-map, then the morphism $E\fun F$ factors through
$L\otimes((\pi_p)_*\cO_{C_p})^*$, and this is not possible since
$\deg(L\otimes((\pi_p)_*\cO_{C_p})^*)>1\geq\deg F$.
If it is not the zero-map, then it is an inclusion because
$(\pi_p)_*\cO_{C_p}$ is torsion-free, and this implies that $\deg F=1$ and 
$F\simeq(\pi_p)_*\cO_{C_p}$.
But this can happen only if the extension we started with is trivial.
\qed

We saw in Proposition \ref{cEL} that there exists a locally-free sheaf
$\cE_L$ on $\bP_L\times C$ such that
$\cE_L|_{\{x\}\times C}\simeq\phi_L(x)$ for every $x\in\bP_L\setminus C$.
To prove Theorem \ref{1}, we introduce a torsion-free sheaf $\cE_{L,1}$ on
$\bP_{L,1}\times C$ which induces the rational map $\phi_{L,1}$.
It is defined by
$$ \cE_{L,1} := \ker \( (\varepsilon_1,1)^*\cE_L \lfun
\bigoplus_{p \in J} \cO_{E_1|_p} \boxtimes (\pi_p)_*\cO_{C_p} \). $$
Note that the map is surjective because, for every node $p$,
$\cE_L|_{\{p\}\times C}$ is isomorphic to $\phi_L(p)$, which surjects onto
$(\pi_p)_*\cO_{C_p}$ by Lemma \ref{emb}.
Moreover, the sheaf $\oplus_{p\in J}\cO_{E_1|_p}\boxtimes(\pi_p)_*\cO_{C_p}$
is supported on $E_1\times C$, and so $\cE_{L,1}$ defines the same map as
$\cE_L$ on $\bP_{L,1}\setminus E_1\simeq\bP_L\setminus J$, i.e.,
$\cE_{L,1}|_{\{x\}\times C}\simeq\cE_L|_{\{\varepsilon_1(x)\}\times C}\simeq
\phi_{L,1}(x)=\phi_L(\varepsilon_1(x))$ for every
$x\in\bP_{L,1}\setminus E_1$.
If $x\in E_1|_p$, we have an exact sequence
$ \cE_{L,1}|_{\{x\} \times C} \fun (\varepsilon_1,1)^*\cE_L|_{\{x\} \times C}
\fun (\pi_p)_*\cO_{C_p} \fun 0, $ which completes to an exact sequence on $C$
$$ 0 \lfun T \lfun \cE_{L,1}|_{\{x\} \times C} \lfun
\cE_L|_{\{p\} \times C} \lfun (\pi_p)_*\cO_{C_p} \lfun 0, $$
where $T$ is the torsion sheaf $\cTor_1^{\bP_{L,1}\times C}
(\cO_{E_1|_p}\boxtimes(\pi_p)_*\cO_{C_p},\cO_{\{x\}\times C})$, that we shall
see to be isomorphic to $(\pi_p)_*\cO_{C_p}$.
Since the kernel of $\cE_L|_{\{p\}\times C}\fun(\pi_p)_*\cO_{C_p}$ is
$L\otimes((\pi_p)_*\cO_{C_p})^*$, $\cE_{L,1}|_{\{x\}\times C}$ is an
extension of $L\otimes((\pi_p)_*\cO_{C_p})^*$ by $(\pi_p)_*\cO_{C_p}$.
Therefore, by Lemma \ref{extss}, $\cE_{L,1}|_{\{x\}\times C}$ is
semi-stable if and only if it does not split as such an extension.
To prove Theorem \ref{1}, we need to show that $\cE_{L,1}|_{E_1\times C}$
induces the rational map $\phi_{L,1}|_{E_1}$, i.e., that
$ \cE_{L,1}|_{\{x\}\times C} \simeq E_x $
for every $x\in E_1$.
This will be proved in Proposition \ref{EL1E1xC}.

\section{Description of $\cE_{L,1}$}\label{cEL1}

The main goal of this section is to prove that $\cE_{L,1}$ induces the
rational map $\phi_{L,1}$, and we start by analyzing $\cE_{L,1}$.
Since $\cE_L$ fits into a short exact sequence
$ 0 \fun \cO_{\bP_L}(1) \fun \cE_L \fun L \fun 0 $
on $\bP_L\times C$, and the image of the composite map
$\varepsilon_1^*\cO_{\bP_L}(1)\inc(\varepsilon_1,1)^*\cE_L\fun
\oplus_{p\in J}\cO_{E_1|_p}\boxtimes(\pi_p)_*\cO_{C_p}$ is
$ \cO_{E_1 \times C}, $
we obtain the following commutative diagram on $\bP_{L,1}\times C$
\begin{eqnarray}\label{diagramE1}
\begin{CD}
& & 0 & & 0 & & 0 \\
& & @VVV @VVV @VVV \\
0 @>>> \cA_1 @>>> \cE_{L,1} @>>> \cB_1 @>>> 0 \\
& & @VVV @VVV @VVV \\
0 @>>> \varepsilon_1^*\cO_{\bP_L}(1) @>>>
(\varepsilon_1,1)^*\cE_L @>>> L @>>> 0 \\
& & @VVV @VVV @VVV \\
0 @>>> \cO_{E_1 \times C} @>>> \ds\bigoplus_{p \in J} \cO_{E_1|_p} \boxtimes
(\pi_p)_*\cO_{C_p} @>>> \ds\bigoplus_{p \in J} \cO_{E_1 \times \{ p \}}@>>>0\\
& & @VVV @VVV @VVV \\
& & 0 & & 0 & & 0
\end{CD},
\end{eqnarray}
where $\cA_1$, $\cE_{L,1}$, and $\cB_1$ are defined by the exactness of
the vertical exact sequences.
In particular, $\cA_1\simeq\pi_{\bP_{L,1}}^*(\varepsilon_1^*
\cO_{\bP_L}(1)\otimes\cO_{\bP_{L,1}}(-E_1))$.

This shows that $\cE_{L,1}$ fits in a short exact sequence
$ 0 \fun \varepsilon_1^*\cO_{\bP_L}(1) \otimes
\cO_{\bP_{L,1}} (-E_1) \fun \cE_{L,1} \fun \cB_1 \fun 0 $
on $\bP_{L,1}\times C$ which restricts to a short exact sequence
$ 0 \fun \cO_{E_1}(1) \fun \cE_{L,1}|_{E_1 \times C}
\fun \cB_1|_{E_1 \times C} \fun 0 $
on $E_1\times C$.
The restriction stays exact because $\cO_{E_1}(1)$ is locally-free, and
the map $\cO_{E_1}(1)\fun\cE_{L,1}|_{E_1\times C}$ is generically
injective.
Therefore, the image of any $\cTor$ sheaf which would appear is $0$.

\begin{rem}
We shall use this fact several times when restricting diagrams or
short exact sequences.
When no comments are made about a sequence staying exact after a
restriction, the reason shall be the same as here, i.e., the first
sheaf is locally-free, and the first map is generically injective.
\end{rem}

\begin{lem}\label{cB1}
For each node $p\in J$, there exists a short exact sequence
$$ 0 \lfun \cO_{E_1|_p}(1) \boxtimes (\pi_p)_*\cO_{C_p} \lfun
\cE_{L,1}|_{E_1|_p \times C} \lfun L \otimes ((\pi_p)_*\cO_{C_p})^* \lfun 0 $$
on $E_1|_p\times C$.
\end{lem}

\pf
If we restrict the diagram (\ref{diagramE1}) to $E_1|_p\times C$, we
obtain
\begin{eqnarray*}
\begin{CD}
& & 0 & & 0 & & 0 \\
& & @VVV @VVV @VVV \\
0 @>>> \cO_{E_1|_p}(1) @>>> \cO_{E_1|_p}(1) \boxtimes (\pi_p)_*\cO_{C_p} @>>>
\cO_{E_1|_p \times \{ p \}} (1) @>>> 0 \\
& & @V \simeq VV @VVV @VVV \\
0 @>>> \cO_{E_1|_p}(1) @>>> \cE_{L,1}|_{E_1|_p \times C} @>>>
\cB_1|_{E_1|_p \times C} @>>> 0 \\
& & @V 0 VV @VVV @VVV \\
0 @>>> \cO_{E_1|_p \times C} @>>> \pi_C^*\cE_L|_{\{p\} \times C} @>>>
L @>>> 0 \\
& & @V \simeq VV @VVV @VVV \\
0 @>>> \cO_{E_1|_p \times C} @>>> \cO_{E_1|_p} \boxtimes (\pi_p)_*\cO_{C_p}
@>>> \cO_{E_1|_p \times \{ p \}} @>>> 0 \\
& & @VVV @VVV @VVV \\
& & 0 & & 0 & & 0
\end{CD}.
\end{eqnarray*}
It follows from the commutativity of the diagram that
$$ \ker(\pi_C^*\cE_L|_{\{p\} \times C} \lfun \pi_C^*((\pi_p)_*\cO_{C_p}))
\simeq \ker(\pi_C^*L \lfun \cO_{E_1|_p \times \{ p \}}) \simeq
\pi_C^*(L \otimes ((\pi_p)_*\cO_{C_p})^*), $$
which implies our statement.
\qed

\begin{pro}\label{EL1E1xC}
The sheaf $\cE_{L,1}$ on $\bP_{L,1}\times C$ induces the rational map
$\phi_{L,1}$.
\end{pro}

\pf
Since we already saw that $\cE_{L,1}$ defines the rational map $\phi_{L,1}$ on
$\bP_{L,1}\setminus E_1$, it suffices to show that, for every node $p$,
$\cE_{L,1}|_{E_1|_p\times C}$ induces the rational map $\phi_{L,1}|_{E_1|_p}$.
Fix a node $p\in J$.

If we pull-back the extension $\cE_L$ to $\bP_{L,1}\times C$, and then
push it forward via the inclusion
$\pi_{\bP_{L,1}}^*(\varepsilon_1^*\cO_{\bP_L}(1))\inc
\varepsilon_1^*\cO_{\bP_L}(1)\boxtimes(\pi_p)_*\cO_{C_p}$
\begin{eqnarray}\label{cE'0}
\begin{CD}
0 @>>> \varepsilon_1^*\cO_{\bP_L}(1) @>>> (\varepsilon_1, 1)^*\cE_L
@>>> L @>>> 0 \\
& & @VVV @VVV || \\
0 @>>> \varepsilon_1^*\cO_{\bP_L}(1) \boxtimes (\pi_p)_*\cO_{C_p} @>>> \cE'_0
@>>> L @>>> 0
\end{CD},
\end{eqnarray}
we obtain an extension $\cE'_0$ which splits when restricted to
$E_1|_p\times C$.
Indeed, $\varepsilon_1^*\cO_{\bP_L}(1)|_{E_1|_p}\simeq\cO_{E_1|_p}$, and
since $\Ext^1_{E_1|_p\times C}(L,(\pi_p)_*\cO_{C_p})\simeq
H^0(E_1|_p,\cO_{E_1|_p})\otimes\Ext^1_C(L,(\pi_p)_*\cO_{C_p})$ (see [Arc04]),
we see that $\cE'_0|_{E_1|_p\times C}$ splits as long as
$\cE'_0|_{\{x\}\times C}$ splits for some $x\in E_1|_p$.
Restricting the diagram (\ref{cE'0}) above to $\{x\}\times C$ for any
$x\in E_1|_p$, we see that $\cE'_0|_{\{x\}\times C}$ is the trivial
extension $\psi_p(\cE_L|_{\{p\}\times C})$.

Therefore, there exists a surjective map
$\cE'_0\fun\cO_{E_1|_p}\boxtimes(\pi_p)_*\cO_{C_p}$: 
Define $\cE'_1$ to be its kernel.

There exists a commutative diagram on $\bP_{L,1}\times C$:
\begin{eqnarray*}
\begin{CD}
& & 0 & & 0 \\
& & @VVV @VVV \\
0 @>>> \cA'_1 @>>> \cE'_1 @>>> L @>>> 0 \\
& & @VVV @VVV || \\
0 @>>> \varepsilon_1^*\cO_{\bP_L}(1) \boxtimes (\pi_p)_*\cO_{C_p} @>>> \cE'_0
@>>> L @>>> 0 \\
& & @VVV @VVV \\
& & \cO_{E_1|_p} \boxtimes (\pi_p)_*\cO_{C_p} & = &
\cO_{E_1|_p} \boxtimes (\pi_p)_*\cO_{C_p} \\
& & @VVV @VVV \\
& & 0 & & 0
\end{CD},
\end{eqnarray*}
with $\cA'_1\simeq(\varepsilon_1^*\cO_{\bP_L}(1)\otimes
\cO_{\bP_{L,1}}(- E_1|_p))\boxtimes(\pi_p)_*\cO_{C_p}$.
Moreover, if we restrict $i_1\colon\cE_{L,1}\inc(\varepsilon_1,1)^*\cE_L$ and
$i'_1\colon\cE'_1\inc\cE'_0$ to $E_1|_p\times C$, we obtain the following
commutative diagram, where the first row is the exact sequence described in
Lemma \ref{cB1}.
\begin{eqnarray*}
\begin{CD}
0 @>>> \cO_{E_1|_p}(1) \boxtimes (\pi_p)_*\cO_{C_p} @>>>
\cE_{L,1}|_{E_1|_p\times C} @> i_1|_{E_1|_p\times C} >>
L \otimes ((\pi_p)_*\cO_{C_p})^* @>>> 0 \\
& & || & & @VVV @VVV \\
0 @>>> \cO_{E_1|_p}(1) \boxtimes (\pi_p)_*\cO_{C_p} @>>> 
\cE'_1|_{E_1|_p\times C} @> i'_1|_{E_1|_p\times C} >> L @>>> 0 \\
\end{CD}.
\end{eqnarray*}
This shows that $\cE_{L,1}|_{E_1|_p\times C}$ is the pull-back of
$\cE'_1|_{E_1|_p\times C}$ via the inclusion
$L\otimes((\pi_p)_*\cO_{C_p})\inc L$ pulled-back from $C$ to $E_1|_p\times C$.

This is a summary of the steps we took in the construction of
$\cE'_1|_{E_1|_p\times C}$:
\begin{eqnarray*}
\begin{matrix}
\cE_L & \in & \Ext^1_{\bP_L \times C}(L, \cO_{\bP_L}(1)) \\
& & \downarrow \\
(\varepsilon_1,1)^*\cE_L & \in & \Ext^1_{\bP_{L,1} \times C}(L,
\varepsilon_1^*\cO_{\bP_L}(1)) \\
& & \downarrow \\
\cE'_0 & \in & \Ext^1_{\bP_{L,1} \times C}(L, \varepsilon_1^*
\cO_{\bP_L}(1) \boxtimes (\pi_p)_*\cO_{C_p}) \\
& & \uparrow \\
\cE'_1 & \in & \Ext^1_{\bP_{L,1} \times C}(L, (\varepsilon_1^* \cO_{\bP_L}(1)
\otimes \cO_{\bP_{L,1}}(- E_1|_p)) \boxtimes (\pi_p)_*\cO_{C_p})\\
& & \downarrow \\
\cE'_1|_{E_1|_p \times C} & \in & \Ext^1_{E_1|_p \times C}(L, \cO_{E_1|_p}(1)
\boxtimes (\pi_p)_*\cO_{C_p})
\end{matrix}
\end{eqnarray*}

Using the natural isomorphisms
$ \Ext^1_{Y \times C} (L, F \boxtimes G)
\simeq H^0(Y, F) \otimes \Ext^1_C(L, G) $
(see [Arc04]), we can understand what extension $\cE'_1|_{E_1|_p\times C}$ is
by tracking the corresponding elements in these spaces.
Let $v_0,\dots,v_n$ be a basis of $\Ext^1_C(L,\cO_C)$, with
$\Span\{v_0\}=\langle p\rangle$, and let $v_0^*,\dots,v_n^*$ be the
corresponding dual basis in
$\Ext^1_C(L,\cO_C)^*\simeq H^0(\bP_L,\cO_{\bP_L}(1))$.
Then $\cE_L$ corresponds to the element
$ \sum_{i=0}^n v_i^*\otimes v_i \in
H^0(\bP_L, \cO_{\bP_L}(1)) \otimes \Ext^1_C(L, \cO_C), $
and
$ \cE'_1|_{E_1|_p \times C} $
corresponds to the element
$ \sum_{i=1}^n \psi_p(v_i)^* \otimes \psi_p(v_i) \in 
H^0(E_1|_p, \cO_{E_1|_p}(1)) \otimes \Ext^1_C(L, (\pi_p)_*\cO_{C_p}). $
Therefore, since $\cE_{L,1}|_{E_1|_p\times C}$ is the pull-back
of $\cE'_1|_{E_1|_p\times C}$ via $L\otimes((\pi_p)_*\cO_{C_p})^*\inc L$,
$\cE_{L,1}|_{E_1|_p\times C}$ corresponds to $\psi_{L_p}$ itself.
This proves that, for any $a\in\Ext^1_C(L,(\pi_p)_*\cO_{C_p})$, $a\neq0$,
$ \cE_{L,1}|_{\{[a]\}\times C} $ is $ \psi_{L_p}(a) $
as extensions of $L\otimes((\pi_p)_*\cO_{C_p})^*$ by $(\pi_p)_*\cO_{C_p}$.
\qed

\section{The second blow-up}

We now blow-up $\bP_{L,1}$ along the lines $L_p\su E_1|_p$ ($p\in J$).
Let
$$ \bP_{L,2} := \BL_{\coprod_{p \in J} L_p} \bP_{L,1}
\overset{\varepsilon_2}{\lfun} \bP_{L,1}\overset{\varepsilon_1}{\lfun}\bP_L,$$
and let $E_2\su\bP_{L,2}$ be the exceptional divisor, which is the disjoint
union of projective bundles $E_{2,p}\fun L_p$ ($p\in J$).

\begin{thm}\label{2}
$(a)$ The composition
$ \phi_{L,1} \circ \varepsilon_2 \colon \bP_{L,2} \lfun \oSU $
extends to a rational map $\phi_{L,2}$ with the following property.
For each $l\in L_p$, the rational map
$ \phi_{L,2}|_{E_2|_l} \colon E_2|_l \lfun \oSU $
is the projectivization of a linear homomorphism
$\cN_{L_p/\bP_{L,1}}|_l\fun H'_p$, where $H'_p$ is the closure of the locus of
vector bundles of determinant $L$ in
$ \eExt^1_C(L \otimes ((\pi_p)_*\cO_{C_p})^*, (\pi_p)_*\cO_{C_p}). $
This linear homomorphism is an isomorphism if $l\ne p_1,p_2$, and it
maps $\cN_{L_p/\bP_{L,1}}|_{p_i}$ $(i=1,2)$ surjectively onto the hyperplane
$\emph{Im}\psi_{L_p}\su H'_p$.

$(b)$ The indeterminancy locus of
$\phi_{L,2}\colon\bP_{L,2}\fun\oSU$ is the strict transform $\tC_2$
of $\tC_1$.
\end{thm}

\begin{cor}\label{imphiL2}
The image $\phi_{L,2}(\bP_{L,2}\setminus\tC_2)$ of $\phi_{L,2}$ in
$\oSU$ is given by
$$ \phi_L(\bP_L \setminus C) \cup \bigcup_{p \in J} \bP(H'_p). $$
\end{cor}

\begin{rem}
Here $\bP(H'_p)$ actually stands for its image in $\oSU$ via the natural
`forgetful' map, which is a morphism by Lemma \ref{extss}.
We shall prove that this morphism is injective if $g>\deg L$.
\end{rem}

Note that the strict transform $\tC_2$ of $\tC_1$ is isomorphic to $\tC_1$ and
$N$, and, for each node $p$, it intersects $E_{2,p}$ at two points $\tp_1$ and
$\tp_2$ lying over $p_1$ and $p_2$, respectively.

The first step in the proof of Theorem \ref{2} is the analysis of the
exceptional divisor $E_2$.
For each node $p$, $E_{2,p}$ is canonically isomorphic to the projective
bundle $\bP(\cN_{L_p/\bP_{L,1}})$ over $L_p$.
Since $\cN_{L_p/\bP_{L,1}}$ is the normal bundle to $L_p$ in $\bP_{L,1}$, it
contains the normal bundle to $L_p$ in $E_1$, and we obtain short exact
sequences
$ 0 \fun \cN_{L_p/E_1} \fun \cN_{L_p/\bP_{L,1}} \fun
\cN_{E_1/\bP_{L,1}}|_{L_p} \fun 0 $
of vector bundles on $L_p$.

\begin{lem}\label{Nsplit}
For each node $p$, the sequence
$ 0 \fun \cN_{L_p/E_1} \fun \cN_{L_p/\bP_{L,1}} \fun
\cN_{E_1/\bP_{L,1}}|_{L_p} \fun 0 $
splits.
If $N=\dim\bP_L$, then
$ \cN_{L_p/E_1} \simeq \cO_{L_p}(1)^{\oplus N-2}, $
and
$ \cN_{E_1/\bP_{L,1}}|_{L_p} \simeq \cO_{L_p}(-1). $
Moreover, $\bP(\cN_{E_1/\bP_{L,1}}|_{L_p})$ maps isomorphically to
$\widetilde{T_pC}\cap E_2$ via $\varepsilon_2\circ\varepsilon_1$, where
$\widetilde{T_pC}$ is the strict transform of $T_pC$ in $\bP_{L,2}$.
\end{lem}

We shall denote $\bP(\cN_{E_1/\bP_{L,1}}|_{L_p})$ by $L_{2,p}$.
It is isomorphic to $L_p$ via $\varepsilon_2|_{L_{2,p}}$, and it
corresponds to a section of $E_{2,p}\fun L_p$.

\pf
The short exact sequence
$ 0 \fun \cI_{L_p}/\cI_{L_p}^2 \fun \Omega_{E_1}|_{L_p} \fun
\Omega_{L_p} \fun 0, $
together with the standard short exact sequence for $\Omega_{\bP^n}$ (see
[Har77, II.8.13]), proves that
$ \cN_{L_p/E_1} \simeq \cO_{L_p}(1)^{\oplus N-2}. $
It is a standard fact about blow-ups that
$\cN_{E_1/\bP_{L,1}}|_{L_p}\simeq\cO_{L_p}(-1)$.
Finally, the short exact sequence
$ 0 \fun \cN_{L_p/E_1} \fun \cN_{L_p/\bP_{L,1}} \fun
\cN_{E_1/\bP_{L,1}}|_{L_p} \fun 0 $
splits because every extension of $\cO_{L_p}(-1)$ by
$\cO_{L_p}(1)^{\oplus N-2}$ splits.
Indeed,
$ \Ext^1_{L_p}(\cO_{L_p}(-1),\cO_{L_p}(1)^{\oplus N-2}) $
is isomorphic to
$ \oplus^{N-2} H^1(L_p,\cO_{L_p}(2)) = 0. $

The prove the last statement of the lemma, note that,
for each $l\in L_p$, if we let $X_l$ be the projective line in $\bP_L$
which passes through $p$ and corresponds to $l$, we obtain the following
canonical isomorphisms:
\begin{eqnarray*}
\begin{CD}
\cN_{E_1/\bP_{L,1}}|_l & \ \simeq \ &
\ds\frac{T_l \bP_{L,1}}{T_l E_1} & \ \overset{\simeq}{\lfun} \ &
T_p X_l & \ \simeq \ & \ds\frac{\langle X_l \rangle}{\langle p \rangle}
& \ \simeq \ & \langle l \rangle \\
& & & & & & \cap & & \cap \\
& & & & & & \ds\frac{\Ext^1_C(L, \cO_C)}{\langle p \rangle} & \simeq &
\Ext^1_C(L, (\pi_p)_*\cO_{C_p})
\end{CD},
\end{eqnarray*}
where, if $V$ is a vector space, and $S$ is a linear subspace of $\bP(V)$, we
denote by $\langle S\rangle$ the linear subspace of $V$ corresponding to $S$.
\qed

Before we proceed to the proof of Theorem \ref{2}, it is important to study
the following situation:
Fix a node $p$ of $C$ (throughout this section), let $T_pC$ be the projective
tangent plane to $C$ at $p$, and let $X$ be a projective line in $T_pC$
passing through $p$.
As we saw in Lemma \ref{L1}, such lines are parametrized by $L_p$.
Any such line $X_l$ ($l\in L_p$) intersects $C$ at $p$ (and possibly at
other points, but always a finite number), and there exists a rational
map
$ \phi_L|_{X_l} \colon X_l \fun \oSU, $
which extends uniquely to a morphism $\psi_l$ defined on the whole
$X_l$.
We are interested in finding $\psi_l(p)$.
The points of $L_p$ are in one-to-one correspondence with torsion-free sheaves
$M_l$ ($l\in L_p$) of rank $1$ and degree $2$ containing $(\pi_p)_*\cO_{C_p}$,
and
$$ X_l = \bP(\ker(\Ext^1_C(L, \cO_C) \lfun \Ext^1_C(L, M_l))). $$
Note that, if $l\neq p_1,p_2$, then $M_l$ is a line bundle, and if $l=p_i$
$(i=1,2)$, then $M_{p_i}$ is $(\pi_p)_*\cO_{C_p}(p_i)$.

\begin{lem}\label{limitl}
Let $l\in L_p$, $l\neq p_1,p_2$.
Then $\psi_l(p)$ is the unique $($up to isomorphisms$)$ torsion-free sheaf
$E_l$ which can be written both as an extension
$$ 0 \lfun (\pi_p)_*\cO_{C_p} \lfun E_l \lfun L \otimes ((\pi_p)_*\cO_{C_p})^*
\lfun 0 $$
and an extension
$$ 0 \lfun L \otimes M_l^*\lfun E_l \lfun M_l \lfun 0. $$
In particular, it is locally-free.
\end{lem}

\pf
Since for every $x\in X_l\setminus\{p\}$, $\psi_l(x)$ maps onto $M_l$, the
same is true for $\psi_l(p)$.
Indeed, it cannot surject onto something of smaller degree, or it
would not be semi-stable.
Since $\psi_l(p)$ is in $\oSU$, the kernel of $\psi_l(p)\fun M_l$ must
then be $L\otimes M_l^*$, and we have a short exact sequence
$0\fun L\otimes M_l^*\fun \psi_l(p)\fun M_l\fun0$.

Since $\psi_l(x)$ surjects onto $L$ for every $x\in X_l\setminus\{p\}$,
$\psi_l(p)$ surjects onto some torsion-free sheaf $F\su L$, which must have
$\deg F\geq2$ because $\psi_l(p)$ is semi-stable.
Moreover, $F\neq L$ because in that case $\psi_l(p)$ would be an
extension of $L$ by $\cO_C$, but we know that the limit in $\bP_L$ is
$\cE_L|_{\{p\}\times C}$, which is not semi-stable.
Since $\psi_l(p)$ is an extension of $M_l$ by $L\otimes M_l^*$, and every map
from $M_l$ to $F$ is zero because $M_l\not\simeq F$, the composite map
$L\otimes M_l^*\fun\psi_l(p)\fun F$ is non-zero, and therefore
$L\otimes M_l^*\su F\su L$, which implies that
$F\simeq L\otimes((\pi_p)_*\cO_{C_p})^*$.

Therefore, $\psi_l(p)$ is both an extension of $M_l$ by $L\otimes M_l^*$ and
of $L\otimes((\pi_p)_*\cO_{C_p})^*$ by $(\pi_p)_*\cO_{C_p}$, as claimed.
Any such sheaf is in the kernel of the natural linear homomorphism
$ \Ext^1_C(M_l, L \otimes M_l^*) \fun
\Ext^1_C(M_l, L \otimes ((\pi_p)_*\cO_{C_p})^*) $
which is one-dimensional, being isomorphic to $\Hom_C(M_l,\bC_p)$.
\qed

We shall prove in Proposition \ref{EL2E2xC} that the points of
$L_{2,p}=\widetilde{T_pC}\cap E_2\simeq L_p$ map to these vector bundles
$E_l$.
The following lemma describes their geometry.

\begin{lem}\label{conic}
The torsion-free shaves $E_l$ $(l\in L_p)$ form a conic in a quadric
$Q$ in
$$ \bP^3 \simeq \bP(\ker(
\eExt^1_C(L \otimes ((\pi_p)_*\cO_{C_p})^*, (\pi_p)_*\cO_{C_p}) \lfun
\eExt^1_C(L \otimes ((\pi_p)_*\cO_{C_p})^*, (\pi_p)_*\cO_{C_p}(p_1+p_2)))). $$
\end{lem}

\pf
For every $l\in L_p$, let
\begin{eqnarray*}
X_{l,1} & := & \bP(\ker(
\Ext^1_C(L \otimes ((\pi_p)_*\cO_{C_p})^*, (\pi_p)_*\cO_{C_p}) \lfun
\Ext^1_C(L \otimes ((\pi_p)_*\cO_{C_p})^*, M_l))), \\
X_{l,2} & := & \bP(\ker(
\Ext^1_C(L \otimes ((\pi_p)_*\cO_{C_p})^*, (\pi_p)_*\cO_{C_p}) \lfun
\Ext^1_C(L \otimes M_l^*, (\pi_p)_*\cO_{C_p}))).
\end{eqnarray*}
For each $l\in L_p$, the lines $X_{l,1}$ and $X_{l',2}$ span the
plane\footnote{For the dimension of the kernels involved in this proof, see
[Arc04].}
$$ \bP(\ker(\Ext^1_C(L \otimes ((\pi_p)_*\cO_{C_p})^*, (\pi_p)_*\cO_{C_p})
\lfun \Ext^1_C(L \otimes M_{l'}^*, M_l))), $$
and the union of all these lines is a quadric $Q$.
We shall show in Lemma \ref{hyper} that
$H'_p:=\overline{\{\det E\simeq L}\}$ is a hyperplane in
$\Ext^1_C(L\otimes((\pi_p)_*\cO_{C_p})^*,(\pi_p)_*\cO_{C_p})$.
Similarly, $\bP(H'_p)\cap\bP^3$ is a hyperplane in this $\bP^3$, and the
intersection of $\bP(H'_p)$ with $Q$ is the conic in the lemma, since we
know that each $E_l$ is contained in both.
\qed

To prove Theorem \ref{2}, we shall first construct a torsion-free sheaf
$\cE_{L,2}$ on $\bP_{L,2}\times C$, and then show that it induces the rational
map $\phi_{L,2}$.
We can construct $\cE_{L,2}$, starting with the torsion-free sheaf $\cE_{L,1}$
corresponding to the rational map $\phi_{L,1}$, as follows\footnote{For the
existence of the map in the definition of $\cE_{L,2}$, for a proof of its
surjectivity, and for a more in depth analysis of the sheaf, see Section
\ref{cEL2}.}
$$ \cE_{L,2} := \ker \( (\varepsilon_2,1)^*\cE_{L,1} \lfun \bigoplus_{p \in J}
\varepsilon_2^*\cO_{L_p}(1) \boxtimes (\pi_p)_*\cO_{C_p} \). $$
Moreover, the sheaf
$\oplus_{p\in J}\varepsilon_2^*\cO_{L_p}(1)\boxtimes(\pi_p)_*\cO_{C_p}$ is
supported on $E_2\times C$, and so $\cE_{L,2}$ defines the same map as
$\cE_{L,1}$ on
$\bP_{L,2}\setminus E_2\simeq\bP_{L,1}\setminus\coprod_{p\in J}L_p$.
The situation is very similar to the one in the first blow-up, and it
is easy to see that, for every node $p$, every $l\in L_p$ and every
$x\in E_2|_l$, we have an exact sequence
$$ 0 \lfun (\pi_p)_*\cO_{C_p} \lfun \cE_{L,2}|_{\{x\} \times C} \lfun
\cE_{L,1}|_{\{l\} \times C} \lfun (\pi_p)_*\cO_{C_p} \lfun 0. $$
Since $\cE_{L,1}|_{\{l\}\times C}\simeq(\pi_p)_*\cO_{C_p}\oplus
(L\otimes((\pi_p)_*\cO_{C_p})^*)$ for every $l\in L_p$, we obtain that, for
every $x\in E_{2,p}$, $\cE_{L,2}|_{\{x\}\times C}$ is an extension of the
following type:
$$ 0 \lfun (\pi_p)_*\cO_{C_p} \lfun \cE_{L,2}|_{\{x\} \times C} \lfun L
\otimes ((\pi_p)_*\cO_{C_p})^* \lfun 0. $$
We shall prove in Proposition \ref{EL2E2xC} that the sheaf $\cE_{L,2}$
induces the rational map $\phi_{L,2}$.
In particular, for every $l\in L_p$, the restriction of $\phi_{L,2}$ to
$E_2|_l$ is a rational map
$$ E_2|_l \lfun \bP(\Ext^1_C(L \otimes ((\pi_p)_*\cO_{C_p})^*,
(\pi_p)_*\cO_{C_p})), $$
that we want to prove to be linear, and to be a morphism for
$l\ne p_1,p_2$.

\begin{lem}\label{linearity}
For every node $p$, and every $l\in L_p$, the rational map
$$ \phi_{L,2}|_{E_2|_l} \colon E_2|_l \lfun
\bP(\eExt^1_C(L \otimes ((\pi_p)_*\cO_{C_p})^*, (\pi_p)_*\cO_{C_p})) $$
is linear.
\end{lem}

\pf
We give a direct proof of this lemma, but it also follows from the fact, that
we shall prove in Lemma \ref{cEL2E2}, that
$\cE_{L,2}|_{E_2|_l\times C}\in\Ext^1_{E_2|_l\times C}
(L\otimes((\pi_p)_*\cO_{C_p})^*,\cO_{E_2|_l}(1)\boxtimes(\pi_p)_*\cO_{C_p})$
(see [Arc04]).

From the proof of Lemma \ref{Nsplit}, it is clear that there exists a
commutative diagram:
\begin{eqnarray*}
\begin{CD}
E_2|_l @> \phi_{L,2}|_{E_2|_l} >>
\bP(\Ext^1_C(L \otimes ((\pi_p)_*\cO_{C_p})^*, (\pi_p)_*\cO_{C_p})) \\
\bigcup & & \bigcup \\
\bP(\cN_{L_p/E_1}|_l) @> \simeq >> \bP(\Im \psi_{L_p})
\end{CD}.
\end{eqnarray*}
Since the morphism $\bP(\cN_{L_p/E_1}|_l)\fun\bP(\Im \psi_{L_p})$ is a
linear isomorphism, the map itself is linear.
\qed

Let us now show that $H'_p$ is a hyperplane in
$\Ext^1_C(L\otimes((\pi_p)_*\cO_{C_p})^*,(\pi_p)_*\cO_{C_p})$.

\begin{lem}\label{hyper}
The closure $H'_p$ of the locus
$ \{ E \in \eExt^1_C(L \otimes ((\pi_p)_*\cO_{C_p})^*, (\pi_p)_*\cO_{C_p})
\mid \det E \simeq L \} $
in $\eExt^1_C(L\otimes((\pi_p)_*\cO_{C_p})^*,(\pi_p)_*\cO_{C_p})$ is a vector
subspace of codimension $1$.
\end{lem}

\pf
It is enough to show that the closure of
$$ \bP \left( \{ E \in \Ext^1_C(L \otimes ((\pi_p)_*\cO_{C_p})^*,
(\pi_p)_*\cO_{C_p}) \mid \det E \simeq L \} \right) $$
in $\bP(\Ext^1_C(L\otimes((\pi_p)_*\cO_{C_p})^*,(\pi_p)_*\cO_{C_p}))$ is a
linear hyperplane.
Let $E$ be a vector bundle in
$\Ext^1_C(L\otimes((\pi_p)_*\cO_{C_p})^*,(\pi_p)_*\cO_{C_p})$:
What are the possible values for $\det E$?

Let $l\in L_p\setminus\{p_1,p_2\}$, and consider the lines $X_{l,1}$ and
$X_{l,2}$ that we defined in Lemma \ref{conic}.
All vector bundles $E$ in $X_{l,1}$ are of the form
$ 0 \fun L \otimes M_{l'}^*\fun E \fun M_l \fun 0 $
for $l'\in L_p\setminus\{p_1,p_2\}$, and their determinant is of the form
$L\otimes M_l\otimes M_{l'}^*$.
Similarly, all of the vector bundles $E$ in $X_{l,2}$ are of the form
$ 0 \fun L \otimes M_l^*\fun E \fun M_{l'} \fun 0, $
and so their determinant is of the same form as above.
Consider the rational map
$$ \det: \bP(\Ext^1_C(L \otimes ((\pi_p)_*\cO_{C_p})^*, (\pi_p)_*\cO_{C_p}))
\lfun \overline{\{ L' \in \Pic^d(C) \mid \pi_p^*L' \simeq \pi_p^*L \}} \simeq
\bP^1. $$
It is defined on the locus of locally-free sheaves, and it extends to the
locus of the extensions which are not push-forwards of extensions from $C_p$
(see [Bho92]).
Since it is an isomorphism on each line $X_{l,i}$ with $l\neq p_1,p_2$
and $i\in\{1,2\}$, it is a surjective linear map.
\qed

\begin{pro}\label{linearity2}
For every node $p$, and every $l\in L_p$, $l\ne p_1,p_2$, the rational map
$$ E_2|_l \fun \bP(\eExt^1_C(L \otimes ((\pi_p)_*\cO_{C_p})^*,
(\pi_p)_*\cO_{C_p})) $$
is an isomorphism onto its image $\bP(H'_p)$.
If $l=p_i$ $(i=1,2)$, then it maps $E_2|_l$ onto $\bP(\emph{Im}\psi_{L_p})$.
In particular, if $l\ne p_1,p_2$, then $\phi_2|_{E_2|_l}$ is a morphism.
\end{pro}

\pf
We already saw in Lemma \ref{linearity} that the map is linear for
every $l\in L_p$.
Let $x_l$ be the point of intersection between the strict transform of the
projective line $X_l$ with $E_2$.
We shall prove in Proposition \ref{EL2E2xC} that, if $l\neq p_1,p_2$,
$\cE_{L,2}|_{\{x_l\}\times C}$ is isomorphic to the vector bundle $E_l$
of Lemma \ref{limitl}.
Therefore, in this case, the image of $E_2|_l$ contains $\bP(\Im\psi_{L_p})$
and $E_l$.
Since $E_l\not\in\bP(\Im\psi_{L_p})$, the image of $E_2|_l$ is the hyperplane
$H'_p$.

If $l=p_i$ ($i=1,2$), then the image is just $\bP(\Im\psi_{L_p})$.
Indeed, we already know that the map cannot be defined everywhere on
$E_2|_{p_i}$ ($i=1,2$) because it contains a point on the strict
transform $\tC_2$ of $C$ which is contained in the locus of
indeterminancy of $\phi_{L,2}$.
Therefore, it cannot be an isomorphism.
Being a linear map, its image is contained in a hyperplane, which has
to be $\bP(\Im\psi_{L_p})$.
\qed

Theorem \ref{2} will now follow from Proposition \ref{EL2E2xC}.

\section{Description of $\cE_{L,2}$}\label{cEL2}

Since, for every node $p$, $\cE_{L,1}|_{L_p\times C}$ splits as
$ (\cO_{L_p}(1) \boxtimes (\pi_p)_*\cO_{C_p}) \oplus
(L \otimes ((\pi_p)_*\cO_{C_p})^*), $
the map
$ \cE_{L,1} \fun \oplus_{p \in J} \cO_{L_p}(1) \boxtimes (\pi_p)_*\cO_{C_p} $
which appears in the definition of $\cE_{L,2}$ is surjective.

Since the image of the composite map
$$ (\varepsilon_2,1)^*\cA_1 \lfun (\varepsilon_2,1)^*\cE_{L,1} \lfun
\bigoplus_{p \in J} \varepsilon_2^*\cO_{L_p}(1) \boxtimes (\pi_p)_*\cO_{C_p}$$
is $\pi_{E_2}^*\varepsilon_2^*\cO_{L_p}(1)$, we obtain the following
commutative diagram on $\bP_{L,2}\times C$:
\begin{eqnarray}\label{diagramE2}
\begin{CD}
& & 0 & & 0 & & 0 \\
& & @VVV @VVV @VVV \\
0 @>>> \cA_2 @>>> \cE_{L,2} @>>> \cB_2 @>>> 0 \\
& & @VVV @VVV @VVV \\
0 @>>> (\varepsilon_2,1)^*\cA_1 @>>> (\varepsilon_2,1)^*\cE_{L,1}
@>>> (\varepsilon_2,1)^*\cB_1 @>>> 0 \\
& & @VVV @VVV @VVV \\
0 @>>> \ds\bigoplus_{p \in J} F_p @>>> \ds\bigoplus_{p \in J} F_p \boxtimes
(\pi_p)_*\cO_{C_p} @>>> \ds\bigoplus_{p \in J} F_p \boxtimes \bC_p @>>> 0 \\
& & @VVV @VVV @VVV \\
& & 0 & & 0 & & 0
\end{CD},
\end{eqnarray}
where $\cA_2$, $\cE_{L,2}$, and $\cB_2$ are defined by the vertical exact
sequences, and we denoted the locally-free sheaf $\varepsilon_2^*\cO_{L_p}(1)$
on $E_{2,p}$ by $F_p$ to simplify the notation.

We want to show that, for every node $p$, and every $l\in L_p$,
$l\ne p_1,p_2$, $\cE_{L,2}|_{E_2|_l\times C}$ is the universal bundle
associated to $H'_p$ when we identify $E_2|_l$ with $H'_p$.
Let us start with a lemma.

\begin{lem}\label{cEL2E2}
For every node $p$, there exists a short exact sequence
$$ 0 \lfun (\varepsilon_2^*\cO_{L_p}(1) \otimes
\cO_{\bP_{L,2}}(-E_2)|_{E_2}) \boxtimes (\pi_p)_*\cO_{C_p} \lfun
\cE_{L,2}|_{E_{2,p} \times C} \lfun L \otimes ((\pi_p)_*\cO_{C_p})^* \lfun 0$$
on $E_{2,p}\times C$.
Moreover, for each $l\in L_p$, there exists a short exact sequence
$$ 0 \lfun \cO_{E_2|_l}(1) \boxtimes (\pi_p)_*\cO_{C_p} \lfun
\cE_{L,2}|_{E_2|_l \times C} \lfun L \otimes ((\pi_p)_*\cO_{C_p})^* \lfun 0 $$
on $E_2|_l\times C$.
\end{lem}

\pf
The restriction of diagram (\ref{diagramE2}) to $E_{2,p}\times C$ is
\begin{eqnarray*}
\begin{CD}
& & 0 & & 0 & & 0 \\
& & @VVV @VVV @VVV \\
0 @>>> F_p(-E_2) @>>> F_p(-E_2) \boxtimes (\pi_p)_*\cO_{C_p} @>>> F_p(-E_2)
\boxtimes \bC_p @>>> 0 \\
& & @V \simeq VV @VVV @VVV \\
0 @>>> F_p(-E_2) @>>> \cE_{L,2}|_{E_{2,p} \times C} @>>>
\cB_2|_{E_{2,p}\times C} @>>> 0 \\
& & @V 0 VV @VVV @VVV \\
0 @>>> F_p @>>> (\varepsilon_2,1)^*\cE_{L,1}|_{L_p \times C} @>>>
(\varepsilon_2,1)^*\cB_1|_{L_p \times C} @>>> 0 \\
& & @V \simeq VV @VVV @VVV \\
0 @>>> F_p @>>> F_p \boxtimes (\pi_p)_*\cO_{C_p} @>>> F_p \boxtimes \bC_p
@>>> 0 \\
& & @VVV @VVV @VVV \\
& & 0 & & 0 & & 0
\end{CD},
\end{eqnarray*}
where $F_p$ is $\varepsilon_2^*\cO_{L_p}(1)$ as above, and $F_p(-E_2)$ is
$\varepsilon_2^*\cO_{L_p}(1)\otimes\cO_{E_2}(-E_2)$.

The first statement of the lemma follows directly from the diagram by looking
at the middle column and observing that the kernel of the map
$(\varepsilon_2,1)^*\cE_{L,1}|_{L_p\times C}\fun
\varepsilon_2^*\cO_{L_p}(1)\boxtimes(\pi_p)_*\cO_{C_p}$ on $E_{2,p}\times C$
is $\pi_C^*(L\otimes((\pi_p)_*\cO_{C_p})^*)$.

For the second statement, the diagram shows that $\cB_2|_{E_{2,p}\times C}$
has torsion.
Also, for every $l\in L_p$, $\cB_2|_{E_2|_l\times C}$ has torsion, and 
$ \cB_2|_{E_2|_l \times C} / \Tors $ is isomorphic to
$\ker ((\varepsilon_2,1)^*\cB_1|_{\{ l \} \times C} \lfun 
\cO_{E_2|_l \times \{ p \}}). $
By the way we defined this last map, it is clear that this kernel is
just the pull-back via $(\varepsilon_2,1)$ of $\cB_1|_{\{l\}\times C}$
modulo torsion, which is $\pi_C^*(L\otimes((\pi_p)_*\cO_{C_p})^*)$.
From the diagram, it is clear that the kernel of the map
$\cE_{L,2}|_{E_2|_l \times C} \lfun \pi_C^*(L \otimes ((\pi_p)_*\cO_{C_p})^*)$
is the same as the kernel of the map
$ \cE_{L,2}|_{E_2|_l \times C} \lfun (\varepsilon_2,1)^*
\cE_{L,1}|_{\{ l \} \times C}, $
which is 
$ \cO_{E_2|_l}(1) \boxtimes (\pi_p)_*\cO_{C_p}. $
\qed

The following proposition will conclude the proof of Theorem \ref{2}.

\begin{pro}\label{EL2E2xC}
The sheaf $\cE_{L,2}$ on $\bP_{L,2}\times C$ induces the rational map
$\phi_{L,2}$.
\end{pro}

Let us start with a lemma.

\begin{lem}\label{codim}
If $Y\su\bP_{L,2}$ is a smooth subvariety such that
$$ \emph{codim}(Y, \bP_{L,2}) = \emph{codim}(Y \cap E_2, E_2), $$
then $\cTor_1^{\bP_{L,2}\times C}(\cO_{E_2\times C},\cO_{Y\times C})=0$
and
$\cTor_1^{\bP_{L,2}\times C}(\cO_{E_2\times\{p\}},\cO_{Y\times C})=0$.
\end{lem}

\pf
Consider the short exact sequence
$ 0 \fun \cO_{\bP_{L,2}}(-E_2) \overset{f}{\fun} \cO_{\bP_{L,2}}
\overset{g}{\fun} \cO_{E_2} \fun 0 $
on $\bP_{L,2}$ and its pull-back to $\bP_{L,2}\times C$.
If we tensor it with $\cO_{Y\times C}$, we obtain the exact sequence
$$ 0 \lfun \cTor_1^{\bP_{L,2} \times C}(\cO_{E_2 \times C},
\cO_{Y \times C}) \lfun \cO_{\bP_{L,2}}(-E_2)|_Y
\overset{f|_{Y \times C}}{\lfun} \cO_{Y \times C}
\overset{g|_{Y \times C}}{\lfun} \cO_{(Y \cap E_2) \times C} \lfun 0 $$
on $Y\times C$, where the zero on the left occurs because
$\cO_{\bP_{L,2}\times C}$ is locally-free.

Since the codimension of $Y\cap E_2$ in $Y$ is $1$, $g$ is zero on the
dense open subset $(Y\setminus(Y\cap E_2))\times C$, and $f$ is an
isomorphism on it.
Therefore, $\cTor_1^{\bP_{L,2}\times C}(\cO_{E_2\times C},
\cO_{Y\times C})$ is supported on $(Y\cap E_2)\times C$, and it must be
zero, being a subsheaf of $\cO_{\bP_{L,2}}(-E_2)|_{Y\times C}$, which
is a locally-free sheaf on a bigger dimensional variety.

Consider now the short exact sequence
$ 0 \fun ((\pi_p)_*\cO_{C_p})^* \fun \cO_C \fun \bC_p \fun 0 $
on $C$ and its pull-back
$ 0 \fun \cO_{E_2} \boxtimes ((\pi_p)_*\cO_{C_p})^* \fun \cO_{E_2 \times C}
\fun \cO_{E_2 \times \{p\}} \fun 0 $
to $E_2\times C$.
If we tensor this exact sequence with $\cO_{Y\times C}$ over
$\cO_{\bP_{L,2}\times C}$, we obtain the exact sequence
$$ 0 \lfun T \lfun \cO_{Y \cap E_2} \boxtimes ((\pi_p)_*\cO_{C_p})^* \lfun
\cO_{(Y\cap E_2)\times C} \lfun \cO_{(Y\cap E_2)\times\{p\}} \lfun 0 $$
on $(Y\cap E_2)\times C$, where $T=\cTor_1^{\bP_{L,2} \times C}
(\cO_{E_2 \times \{p\}},\cO_{Y \times C})$ and the zero on the left
is the sheaf $\cTor_1^{\bP_{L,2}\times C}(\cO_{E_2\times C},\cO_{Y\times C})$.
Just as above\footnote{For another proof of $T$ being $0$, note that
the map
$\cO_{Y\cap E_2}\boxtimes((\pi_p)_*\cO_{C_p})^*\fun\cO_{(Y\cap E_2)\times C}$
is injective because it is the pull-back of the injective map
$((\pi_p)_*\cO_{C_p})^*\fun\cO_C$ via the flat morphism
$(Y\cap E_2)\times C\fun C$.},
$\cO_{Y\cap E_2}\boxtimes((\pi_p)_*\cO_{C_p})^*\fun\cO_{(Y\cap E_2)\times C}$
is an isomorphism on the dense open subset
$(Y\cap E_2)\times(C\setminus\{p\})$, whose complement has codimension
$1$, and therefore $\cTor_1^{\bP_{L,2}\times C}(\cO_{E_2\times\{p\}},
\cO_{Y\times C})$ must be zero, being supported on
$(Y\cap E_2)\times\{p\}$ and contained in the torsion-free sheaf
$\cO_{Y\cap E_2}\boxtimes((\pi_p)_*\cO_{C_p})^*$, which is supported on a
bigger dimensional variety.
\qed

{\em Proof $($of Proposition $\ref{EL2E2xC})$.}
It is clear that $\cE_{L,2}$ defines $\phi_{L,2}$ on $\bP_{L,2}\setminus E_2$.
On $E_2$, we shall divide the proof in two part.
We shall first show that, for every node $p$, $\cE_{L,2}$ defines the rational
map $\phi_{L,2}$ on $\tE_1\cap E_{2,p}\simeq\bP(\cN_{L_p/E_1})\su E_2$, where
$\tE_1$ is the strict transform of $E_1$, and then we prove that, if
$l\neq p_1,p_2$, then $\cE_{L,2}|_{\{x_l\}\times C}\simeq E_l$, where $E_l$ is
the vector bundle of Lemma $\ref{limitl}$.

Since $\phi_{L,2}$ agrees with the rational map defined by $\cE_{L,2}$ on a
dense open subset, we have that $\phi_{L,2}(x)$ is
$\cE_{L,2}|_{\{x\}\times C}$ whenever this is semi-stable.
In particular, the proposition will follow from the fact that $\cE_{L,2}$ is
semi-stable for every $x\in E_{2,p}$ except for $x=\tp_i$, $i=1,2$, which are
the only two points on $E_{2,p}$ where we know that $\phi_{L,2}$ cannot be
defined.

To prove that $\cE_{L,2}$ defines the rational map $\phi_{L,2}$ on
$\tE_1\cap E_{2,p}\simeq\bP(\cN_{L_p/E_1})\su E_2$, restrict the commutative
diagram (\ref{diagramE2}) to $\tE_1|_p\times C$ to obtain:
\begin{eqnarray*}
\begin{CD}
& & 0 & & 0 & & 0 \\
& & @VVV @VVV @VVV \\
0 @>>> \cA_2|_{\tE_1|_p \times C} @>>> \cE_{L,2}|_{\tE_1|_p \times C} @>>>
\cB_2|_{\tE_1|_p \times C} @>>> 0 \\
& & @VVV @VVV @VVV \\
0 @>>> \sigma^*\cO_{E_1|_p}(1) @>>> (\sigma,1)^*\cE_{L,1}|_{E_1|_p \times C}
@>>> (\sigma,1)^*\cB_1|_{E_1|_p \times C} @>>> 0 \\
& & @VVV @VVV @VVV \\
0 @>>> \sigma^*\cO_{L_p}(1) @>>> \sigma^*\cO_{L_p}(1) \boxtimes
(\pi_p)_*\cO_{C_p} @>>> \sigma^*\cO_{L_p}(1) \boxtimes \bC_p @>>> 0 \\
& & @VVV @VVV @VVV \\
& & 0 & & 0 & & 0
\end{CD},
\end{eqnarray*}
where $\sigma\colon\tE_1|_p\fun E_1|_p$ is the restriction of $\varepsilon_2$
to $\tE_1|_p$.
The vertical columns are exact because
$\cTor_1^{\bP_{L,2}\times C}(\varepsilon_2^*\cO_{L_p}(1),
\cO_{\tE_1|_p\times C})=\cTor_1^{\bP_{L,2}\times C}
(\varepsilon_2^*\cO_{L_p}(1)\boxtimes\bC_p,\cO_{\tE_1|_p\times C})=0$.
This is true because $\varepsilon_2^*\cO_{L_p}(1)\simeq
\varepsilon_2^*\cO_{\bP_{L,1}}(-E_1)\otimes\cO_{E_{2,p}\times C}$, and since
$\varepsilon_2^*\cO_{\bP_{L,1}}(-E_1)$ is a locally-free sheaf, it is enough
to show that
$\cTor_1^{\bP_{L,2}\times C}(\cO_{E_{2,p}\times C},\cO_{\tE_1|_p\times C})$
and
$\cTor_1^{\bP_{L,2}\times C}(\cO_{E_{2,p}\times\{p\}},\cO_{\tE_1|_p\times C})$
are both $0$, which was proved in Lemma \ref{codim}.

Since $(\sigma,1)^*\cE_{L,1}|_{E_1|_p\times C}$ is an extension of
$\pi_C^*(L\otimes((\pi_p)_*\cO_{C_p})^*)$ by
$\sigma^*\cO_{E_1|_p}(1)\boxtimes(\pi_p)_*\cO_{C_p}$, there exists a
commutative diagram
\begin{eqnarray*}
\begin{CD}
& & 0 & & 0 \\
& & @VVV @VVV \\
0 & \ \lfun \ & \cA'_2 & \ \lfun \ & \cE_{L,2}|_{\tE_1|_p\times C} &
\ \lfun \ & L \otimes ((\pi_p)_*\cO_{C_p})^* & \ \lfun \ & 0 \\
& & @VVV @VVV || \\
0 & \lfun & \sigma^*\cO_{E_1|_p}(1) \boxtimes (\pi_p)_*\cO_{C_p} & \lfun &
(\sigma,1)^*\cE_{L,1}|_{E_1|_p\times C} & \lfun &
L \otimes ((\pi_p)_*\cO_{C_p})^* & \lfun & 0 \\
& & @VVV @VVV \\
& & \sigma^*\cO_{L_p}(1) \boxtimes (\pi_p)_*\cO_{C_p} & = &
\sigma^*\cO_{L_p}(1) \boxtimes (\pi_p)_*\cO_{C_p} \\
& & @VVV @VVV \\
& & 0 & & 0
\end{CD},
\end{eqnarray*}
where $\cA'_2\simeq(\sigma^*\cO_{E_1|_p}(1)\otimes
\cO_{\tE_1|_p}(-(\tE_1|_p\cap E_2)))\boxtimes(\pi_p)_*\cO_{C_p}$.

Using the natural isomorphisms
$ \Ext^1_{Y \times C} (L, F \boxtimes G)
\simeq H^0(Y, F) \otimes \Ext^1_C(L, G) $
(see [Arc04]) as we did in the proof of Proposition \ref{EL1E1xC}, we have the
following diagram\footnote{Note that $\psi_{L_p}(w_1)=\psi_{L_p}(w_2)=0$.}
\begin{eqnarray*}
\begin{matrix}
\cE_{L,1}|_{E_1|_p \times C} & \longleftrightarrow &
\ds\sum_{i=1}^n w_i^*\otimes \psi_{L_p}(w_i) & \in &
H^0(E_1|_p, \cO_{E_1|_p}(1)) \otimes V \\
& & & & \downarrow \\
(\sigma,1)^*\cE_{L,1}|_{E_1|_p \times C} & \longleftrightarrow &
\ds\sum_{i=1}^n w_i^*\otimes \psi_{L_p}(w_i) & \in &
H^0(\tE_1|_p, \sigma^*\cO_{E_1|_p}(1)) \otimes V \\
& & & & \uparrow \\
\cE_{L,2}|_{\tE_1|_p \times C} & \longleftrightarrow &
\ds\sum_{i=3}^n w_i^*\otimes \psi_{L_p}(w_i) & \in & H^0(\tE_1|_p, \cA'_2)
\otimes V \\
& & & & \downarrow \\
\cE_{L,2}|_{(\tE_1|_p \cap E_2|_l) \times C} & \longleftrightarrow &
\ds\sum_{i=3}^n \psi_{L_p}(w_i)^* \otimes \psi_{L_p}(w_i) & \in &
H^0(\tE_1|_p \cap E_2|_l, \cO_{\tE_1|_p \cap E_2|_l}(1)) \otimes V
\end{matrix},
\end{eqnarray*}
where $w_1,\dots,w_n$ is a basis of $\Ext^1_C(L,(\pi_p)_*\cO_{C_p})$ such that
$\Span\{w_1,w_2\}=\ker\psi_{L_p}$, $w_1^*,\dots,w_n^*$ is the corresponding
dual basis of
$\Ext^1_C(L,(\pi_p)_*\cO_{C_p})^*\simeq H^0(E_1|_p,\cO_{E_1|_p}(1))$, and we
denoted $\Ext^1_C(L \otimes ((\pi_p)_*\cO_{C_p})^*, (\pi_p)_*\cO_{C_p})$ by
$V$ to simplify the diagram.

This proves that the torsion-free sheaf
$\cE_{L,2}|_{(\tE_1|_p\cap E_2|_l)\times C}$ corresponds to the inclusion when
we identify the vector space
$ \Ext^1_{(\tE_1|_p\cap E_2|_l)\times C}(L \otimes ((\pi_p)_*\cO_{C_p})^*,
\cO_{\tE_1|_p \cap E_2|_l}(1) \boxtimes (\pi_p)_*\cO_{C_p}) $
with the vector space
$ \Hom(\Im \psi_{L_p}, \Ext^1_C(L \otimes ((\pi_p)_*\cO_{C_p})^*,
(\pi_p)_*\cO_{C_p})). $
In particular, for every $a\in\Im\psi_{L_p}$, $a\neq0$, $[a]\in
\bP(\Im\psi_{L_p})\simeq\bP(\cN_{L_p/E_1|_p}|_l)\simeq\tE_1|_p\cap E_2|_l$,
and
$ \cE_{L,2}|_{\{[a]\} \times C} \simeq a $
as extensions of $L\otimes((\pi_p)_*\cO_{C_p})^*$ by $(\pi_p)_*\cO_{C_p}$.

To prove that, for every $l\in L_p$, $l\neq p_1,p_2$,
$\cE_{L,2}|_{\{x_l\}\times C}\simeq E_l$, and conclude the proof of the
proposition, it suffices to show that $\cE_{L,2}|_{\{x_l\}\times C}$ is
semi-stable.
This is done by tracking the restrictions of $\cE$, $\cE_1$, and $\cE_2$ to
the product of $X_l$ and its strict transforms with the curve $C$.
Restricting the diagrams defining $\cE_1$ and $\cE_2$ to $\tX_l\times C$, we
obtain a short exact sequence
$$ 0 \lfun L \otimes M_l^* \lfun \cE_2|_{\tX_l \times C} \lfun
\cO_{\tX_l}(-1) \boxtimes M_l \lfun 0. $$
Therefore, $\cE_{L,2}|_{\{x_l\}\times C}$ cannot split, being an extension of
$M_l$ by $L\otimes M_l^*$ and an extension of $L\otimes((\pi_p)_*\cO_{C_p})^*$
by $(\pi_p)_*\cO_{C_p}$, and therefore it is semi-stable.
By continuity, it must be isomorphic to $E_l$, that we proved to be
the limit of $\psi_l(x)$ as $x\mapsto p$.
\qed

\section{The third blow-up}\label{third}

To resolve the indeterminancy of $\phi_{L,2}$, we now blow-up
$\bP_{L,2}$ along $\tC_2$.
Let
$$ \bP_{L,3} := \BL_{\tC_2} \bP_{L,2} \overset{\varepsilon_3}{\lfun} 
\bP_{L,2} \overset{\varepsilon_2}{\lfun} \bP_{L,1}
\overset{\varepsilon_2}{\lfun} \bP_L, $$
and let $E_3\su\bP_{L,3}$ be the exceptional divisor.
For each node $p$, let $\tp_1,\tp_2$ be the points in $\tC_2$ which map to
$p_1,p_2\in\tC_1$, respectively.

\begin{thm}\label{3}
The composition
$ \phi_{L,2} \circ \varepsilon_3 \colon \bP_{L,3} \lfun \oSU $
extends to a morphism $\phi_{L,3}$ such that for each $q\in\tC_2$ not lying
above a node of $C$, the restriction of $\phi_{L,3}$ to $E_3|_q$ maps $E_3|_q$
isomorphically onto\footnote{We identify here a point $q$ on $\tC_2$,
$q\ne\tp_1,\tp_2$, with its image $q$ on $C$.}
$\bP(\eExt^1_C(L(-q),\cO_C(q)))$, and for each node $p\in C$, its restriction
to $E_3|_{\tp_i}$ sends $E_3|_{\tp_i}$ isomorphically onto $\bP(H'_p)$
$(i=1,2)$.
\end{thm}

\begin{cor}
The image of $\phi_{L,3}$ in $\oSU$ is given by\footnote{Note
that by $\bP(H'_p)$ $[$resp.\ $\bP(\Ext^1_C(L(-q),\cO_C(q)))]$ we actually
mean its image into $\oSU$ by the morphism described in the remark after
Corollary \ref{imphiL2} $[$resp.\ in Corollary $\ref{qsuoSU}]$.}
$$ \phi_L(\bP_L \setminus C) \cup \bigcup_{p \in J} \bP(H'_p) \cup
\bigcup_{q \not\in J} \bP(\eExt^1_C(L(-q),\cO_C(q))). $$
\end{cor}

We shall prove Theorem \ref{3} in the next section, after we study the
exceptional divisor $E_3$ in this section.
We know that $E_3$ is canonically isomorphic to $\cN_{\tC_2/\bP_{L,2}}$.

Let $q$ be a point of $\tC_2$ not lying above a node of $C$.
Then we have the canonical isomorphisms
$$ \cN_{\tC_2 / \bP_{L,2}}|_q \simeq \frac{T_q \bP_{L,2}}{T_q \tC_2}
\simeq \frac{T_q \bP_L}{T_q C} \simeq
\frac{\ds\frac{\Ext^1_C(L, \cO_C)}{\langle q \rangle}}{T_q C}, $$
and so, to prove that
$ E_3|_q \overset{\simeq}{\lfun} \bP(\Ext^1_C(L(-q),\cO_C(q))), $
it is necessary to prove that
$$ T_q C \simeq \frac{\ker(\Ext^1_C(L, \cO_C) \lfun
\Ext^1_C(L(-q), \cO_C(q)))}{\langle q \rangle}. $$
Since, as in the proof of Lemma \ref{TpC}, the secant line joining two smooth
points $q,q'$ of $C$ is
$ \bP(\ker(\Ext^1_C(L, \cO_C) \lfun \Ext^1_C(L, \cO_C(q + q')))), $
and $\Ext^1_C(L,\cO_C(q + q'))\simeq\Ext^1_C(L(-q'),\cO_C(q))$, this
follows when taking the limit as $q'\fun q$.

Let now $p$ be a node of $C$, and $q=\tp_i$ with $i\in\{1,2\}$.
Then
$$ \cN_{\tC_2 / \bP_{L,2}}|_{\tp_i} \simeq
\frac{T_{\tp_i} \bP_{L,2}}{T_{\tp_i} \tC_2} \simeq T_{\tp_i} E_2. $$
This contains the canonical hyperplane $T_{\tp_i}(E_2|_{p_i})$ which maps
isomorphically to $\Im\psi_{L_p}$.
Indeed, using Lemma \ref{Nsplit} and the fact that
$E_2\simeq\bP(\cN_{L_p/\bP_{L,1}})$,
$$ T_{\tp_i}(E_2|_{p_i}) \simeq
\frac{\cN_{L_p/\bP_{L,1}}|_{p_i}}{\langle \tp_i \rangle} \simeq 
\frac{\cN_{L_p/E_1}|_{p_i} \oplus
\cO_{L_p}(-1)|_{p_i}}{\cO_{L_p}(-1)|_{p_i}} \simeq
\cN_{L_p/E_1}|_{p_i}, $$
that we already saw to be canonically isomorphic to $\Im\psi_{L_p}$.
We shall see in Proposition \ref{EL3Impsi} that the morphism
$\bP(T_{\tp_i}(E_2|_{p_i}))\fun\oSU$ factors through this canonical
isomorphism, i.e., there exists a commutative diagram
\begin{eqnarray}\label{E3tpi}
\begin{CD}
E_3|_{\tp_i} & \ \simeq \ & \bP(T_{\tp_i} E_2) @>>>
\phi_{L,3}(E_3|_{\tp_i}) & \ \su \ & \oSU \\
& & \cup & & \uparrow \\
& & \bP(T_{\tp_i}(E_2|_{p_i})) @> \simeq >> \bP(\Im\psi_{L_p})
\end{CD}.
\end{eqnarray}
We shall then show that the top map factors through an isomorphism
$E_3|_{\tp_i}\overset{\simeq}{\fun}\bP(H')$.

As for the other blow-ups, to prove Theorem \ref{3}, the strategy is to
construct a universal sheaf $\cE_{L,3}$ on $\bP_{L,3}\times C$, and then prove
that $\cE_{L,3}$ induces the correct rational map.
In this case, we also want to prove that $\cE_{L,3}$ induces a
morphism, i.e., that $\cE_{L,3}|_{\{x\}\times C}$ is semi-stable for
every $x\in\bP_{L,3}$.
The definition of $\cE_{L,3}$ is not as evident as in the other two
blow-ups, and we postpone it to the next section.
By construction, $\cE_{L,3}$ shall agree with $\cE_{L,2}$ on
$\bP_{L,3}\setminus E_3=\bP_{L,2}\setminus\tC_2$, and we shall show in
Propositions \ref{phiL3E3qxC} and \ref{EL3E3qxC} that, if $q\in\tC_2$ does not
lie over a node of $C$, then $\cE_{L,3}|_{E_3|_q\times C}$ induces the
isomorphism of $E_3|_q$ with $\bP(\Ext^1_C(L(-q),\cO_C(q)))$ described above.
To prove that this induces a morphism from $E_3|_q$ to $\oSU$, we need
to prove the following result.

\begin{lem}
All non-trivial extensions in $\eExt^1_C(L(-q),\cO_C(q))$ are semi-stable.
\end{lem}

\pf
This proof is identical to the one of Lemma \ref{extss}.
\qed

\begin{cor}\label{qsuoSU}
The natural `forgetful' map $\bP(\eExt^1_C(L(-q),\cO_C(q)))\fun\oSU$ which
sends an extension $0\fun\cO_C(q)\fun E\fun L(-q)\fun0$ to $E$ is a morphism.
\end{cor}

Fix now a node $p\in C$.
From the definition of $\cE_{L,3}$, it will be clear that, as in the
case of the first two blow-ups, $\phi_{L,3}(E_3|_{\tp_i})\su
\bP(\Ext^1_C(L\otimes((\pi_p)_*\cO_{C_p})^*,(\pi_p)_*\cO_{C_p}))$, and
therefore, using diagram (\ref{E3tpi}), we can prove the following linearity
result.

\begin{lem}
For $i=1,2$, the rational map
$$ \phi_{L,3}|_{E_3|_{\tp_i}} \colon E_3|_{\tp_i} \lfun
\bP(\eExt^1_C(L \otimes ((\pi_p)_*\cO_{C_p})^*, (\pi_p)_*\cO_{C_p})) $$
is linear.
\end{lem}

\pf
The proof is the same as the one of Lemma \ref{linearity}.
\qed

For each $i\in\{1,2\}$, since the map is linear, and we know it to send the
hyperplane $\bP(T_{\tp_i}(E_2|_{p_i}))$ isomorphically onto
$\bP(\Im\psi_{L_p})$, to prove that it maps $E_3|_{\tp_i}$ isomorphically onto
$\bP(H'_p)$, it suffices to show that there exists a point $y\in E_3|_{\tp_i}$
which maps to some point $x\in\bP(H'_p)\setminus\bP(\Im\psi_{L_p})$.

For each point $x\in\bP(H'_p)\setminus\bP(\Im\psi_{L_p})$, there exists a
section $s_x$ of $E_{2,p}\fun L_p$ defined as follows:
If $l\neq p_1,p_2$, $s_x(l)$ is the unique point of $E_2|_l$ which maps to $x$
in $\bP(H'_p)$.
This defines a section on $L_p\setminus\{p_1,p_2\}$, which can be completed to
a section of $E_{2,p}\fun L_p$ by taking its closure.
Note that its closure must satisfy $s_x(p_i)=\tp_i$ for $i=1,2$, because
$\tp_i$ is the only point on $E_2|_{p_i}$ which does not map to
$\bP(\Im\psi_{L_p})$, and $x\not\in\bP(\Im\psi_{L_p})$.

We shall prove in Proposition \ref{E3l} that, for every $l\in L_p$,
$l\neq p_1,p_2$, the point $y_{l,i}$ defined as the intersection of the strict
transform of $s_{E_l}(L_p)$ with $E_3|_{\tp_i}$ maps to $E_l$ for $i=1,2$, and
this shall complete the proof of Theorem \ref{3}.

\section{Definition of $\cE_{L,3}$}\label{cEL3}

We shall define $\cE_{L,3}$ as the kernel of a map
$(\varepsilon_3,1)^*\cE_{L,2}\fun
(\varepsilon_3,1)^*(\cA_2|_{\tC_2\times C}\otimes\cF_2)$, with $\cF_2$
a sheaf on $\tC_2\times C$ such that
$\cF_2|_{\{q\}\times C}\simeq\cO_C(q)$ if $q\neq\tp_1,\tp_2$, and
$\cF_2|_{\{\tp_i\}\times C}\simeq(\pi_p)_*\cO_{C_p}$ for $i=1,2$.
The map corresponds to a map
$\cE_L\fun\pi_{\bP_L}^*\cO_{\bP_L}(1)|_{C\times C}\otimes\cF$, where
$\cF:=\cI_\Delta^*$, $\cI_\Delta$ being the ideal sheaf of the diagonal
$\Delta$ in $C\times C$.

Before we define the map, let us study $\cF$ in more detail.

\begin{lem}\label{cF}
There exists a short exact sequence
$$ 0 \lfun \cO_{C \times C} \lfun \cF \lfun \om_\Delta^{-1} \lfun 0. $$
Moreover,
$\cF|_{\{q\}\times C}\simeq\cO_C(q)$ if $q\neq p$, and 
$\cF|_{\{p\}\times C}\simeq(\pi_p)_*\cO_{C_p}$.
\end{lem}

\pf
Starting with the short exact sequence
$ 0 \fun \cI_\Delta \fun \cO_{C \times C} \fun \cO_\Delta \fun 0, $
and applying the functor $\cHom_{C\times C}(-,\cO_{C\times C})$,
we obtain the short exact sequence\footnote{The sequence starts with
$\cHom_{C\times C}(\cO_\Delta,\cO_{C\times C})$ which is zero
because $\cO_{C\times C}$ is torsion-free, and ends with
$\cExt^1_{C\times C}(\cO_{C\times C},\cO_{C\times C})$ which is also
zero (see [Har77, III.6.3]).}
$$ 0 \lfun \cO_{C \times C} \lfun \cF \lfun
\cExt^1_{C \times C}(\cO_\Delta, \cO_{C \times C}) \lfun 0. $$
Moreover, $\cExt^1_{C\times C}(\cO_\Delta,\cO_{C\times C})\simeq
\cExt^1_{C\times C}(\cO_\Delta,\om_{C\times C})\otimes
\om_{C\times C}^{-1}\simeq\om_\Delta\otimes\om_{C\times C}^{-1}\simeq
\om_\Delta^{-1}$ (see [Har77, III.6.7]) since
$\om_\Delta\simeq\cExt^1_{C\times C}(\cO_\Delta,\om_{C\times C})$
(see [Eis95, 21.15]), and
$\om_{C\times C}|_\Delta\simeq\om_\Delta^{\otimes2}$.

Now, for any $q\in C$, restricting the short exact sequence to
$\{q\}\times C$, we obtain a short exact sequence
$ 0 \fun \cO_C \fun \cF|_{\{q\} \times C} \fun \bC_q \fun 0 $
which does not split.

We know that $\cE_L$ is the extension of $\pi_C^*L$ by
$\pi_{\bP_L}^*\cO_{\bP_L}(1)$ which corresponds to the identity in
$\Hom(\Ext^1_C(L,\cO_C),\Ext^1_C(L,\cO_C))$.
Since
$\pi_{\bP_L}^*\cO_{\bP_L}(1)|_{C\times C}\simeq\pi_1^*(L\otimes\om_C)$,
where $\pi_1$ is the first projection $C\times C\fun C$, we obtain the
following short exact sequence on $C\times C$:
$$ 0 \lfun \pi_1^*(L \otimes \om_C) \lfun \cE_L|_{C \times C} \lfun
\pi_2^* L \lfun 0. $$
The map $\pi_1^*(L\otimes\om_C)\inc\pi_1^*(L\otimes\om_C)\otimes\cF$
extends to a map
$\cE_L|_{C\times C}\fun\pi_1^*(L\otimes\om_C)\otimes\cF$ if
$\cE_L|_{C\times C}$ is in the kernel of the natural linear
homomorphism
$$ \Ext^1_{C \times C}(\pi_2^* L, \pi_1^*(L \otimes \om_C)) \lfun
\Ext^1_{C \times C}(\pi_2^* L, \pi_1^*(L \otimes \om_C) \otimes \cF),$$
i.e., if $\cE_L|_{C\times C}$ is in the image of the natural linear
homomorphism
$$ \Hom_{C \times C}(\pi_2^* L, \pi_1^*(L \otimes \om_C) \otimes
\om_\Delta^{-1}) \lfun \Ext^1_{C \times C}(\pi_2^* L,
\pi_1^*(L \otimes \om_C)). $$
Let us prove that this is the case.
Since $\pi_1^*(L\otimes\om_C)\otimes\om_\Delta^{-1}$ is isomorphic to
$L$ on $\Delta\simeq C$, $\Hom_{C \times C}(\pi_2^*L,
\pi_1^*(L\otimes\om_C)\otimes\om_\Delta^{-1})\simeq
H^0(\Delta,\cO_\Delta)\simeq\bC$.
Moreover,
$$ \Ext^1_{C \times C}(\pi_2^* L, \pi_1^* (L \otimes \om_C)) \simeq
\Ext^1_C(L, \cO_C)^* \otimes \Ext^1_C(L, \cO_C), $$
which has the canonical identity element corresponding to
$\cE_L|_{C\times C}$.
The constant section $1$ of $\cO_\Delta$ maps to the identity, and our claim
is proved, i.e., there exists a map
$ \cE_L|_{C \times C} \fun \pi_{\bP_L}^* \cO_{\bP_L}(1)|_{C \times C}
\otimes \cF $
as claimed at the beginning of the section.

This map is surjective because its restriction to $\{q\}\times C$,
$q\neq p$, [resp.\ to $\{p\}\times C$] is the surjective map
$\cE_L|_{\{q\}\times C}\fun\cO_C(q)$ [resp.\
$\cE_L|_{\{p\}\times C}\fun(\pi_p)_*\cO_{C_p}$] which makes
$\cE_L|_{\{q\}\times C}$ [resp.\ $\cE_L|_{\{p\}\times C}$] not
semi-stable.

There exists a commutative diagram
\begin{eqnarray*}
\begin{CD}
\cE_{L,1}|_{\tC_1 \times C} @>>> \cA_1|_{\tC_1 \times C} \otimes \cF_1
@>>> 0 \\
@VVV @VVV \\
(\sigma,1)^*(\cE_L|_{C \times C}) @> g >>
\sigma^* \cO_{\bP_L}(1)|_{\tC_1 \times C} \otimes (\sigma,1)^*\cF @>>>0
\end{CD},
\end{eqnarray*}
where $\sigma:\tC_1\fun C$ is the restriction of $\varepsilon_1$ to
$\tC_1$, and $\cF_1$ is defined by the first row being exact.
Since the restriction of the short exact sequence defining $\cE_{L,1}$ to
$\tC_1\times C$ stays exact (see Lemma \ref{codim}) the cokernel of the
vertical map on the left is $\cO_{\{p_1,p_2\}}\boxtimes(\pi_p)_*\cO_{C_p}$.
Since the fiber of $\ker g$ at $\{p_i\}\times C$ has degree $2$ for
$i=1,2$, it maps to zero into $(\pi_p)_*\cO_{C_p}$, and we obtain the
following commutative diagram on $\tC_1\times C$:
\begin{eqnarray*}
\begin{CD}
& & & & 0 & & 0 \\
& & & & @VVV @VVV \\
0 @>>> \ker g @>>> \cE_{L,1}|_{\tC_1 \times C} @>>> \cA_1|_{\tC_1
\times C} \otimes \cF_1 @>>> 0 \\
& & || & & @VVV @VVV \\
0 @>>> \ker g @>>> (\sigma,1)^*(\cE_L|_{C \times C}) @> g >> \sigma^*
\cO_{\bP_L}(1)|_{\tC_1} \otimes (\sigma,1)^*\cF @>>> 0 \\
& & & & @VVV @VVV \\
& & & & \cO_{\{p_1,p_2\}} \boxtimes (\pi_p)_*\cO_{C_p} & = & \cO_{\{p_1,p_2\}}
\boxtimes (\pi_p)_* \cO_{C_p} \\
& & & & @VVV @VVV \\
& & & & 0 & & 0
\end{CD}.
\end{eqnarray*}
If we restrict to $\{q\}\times C$, for $q\in\tC_1$, $q\neq p_1,p_2$,
then $\cF_1|_{\{q\}\times C}\simeq\cF|_{\{q\}\times C}\simeq\cO_C(q)$.
Let $i\in\{1,2\}$.
If we restrict the right column to $\{p_i\}\times C$, we obtain
$$ 0 \lfun T \lfun \cF_1|_{\{p_i\} \times C} \lfun \cF|_{\{p\} \times C}
\overset{\simeq}{\lfun} (\pi_p)_* \cO_{C_p} \lfun 0, $$
where
$T=\cTor_1^{\tC_1\times C}(\cO_{\{p_1,p_2\}}\boxtimes(\pi_p)_*\cO_{C_p},
\cO_{\{p_i\}\times C}$).

To calculate this sheaf, consider
$0\fun\cO_{\tC_1}(-p_i)\fun\cO_{\tC_1}\fun\cO_{\{p_i\}}\fun0$ on
$\tC_1$ and its pull-back to $\tC_1\times C$.
We want to tensor it with $\cO_{\{p_1,p_2\}}\boxtimes(\pi_p)_*\cO_{C_p}$, and
we do it in two steps.
We first tensor it with $\pi_C^*((\pi_p)_*\cO_{C_p})$ to obtain
$$ 0 \lfun \cO_{\tC_1}(-p_i) \boxtimes (\pi_p)_* \cO_{C_p} \lfun
(\pi_p)_* \cO_{C_p} \lfun \cO_{\{p_i\}} \boxtimes (\pi_p)_* \cO_{C_p}\lfun0,$$
where the map
$\cO_{\tC_1}(-p_i)\boxtimes(\pi_p)_* \cO_{C_p}\fun(\pi_p)_*\cO_{C_p}$ is
injective because it is an isomorphism on the dense open subset
$(\tC_1\setminus\{p_i\})\times C$ whose complement has codimension $1$,
and therefore the image of any torsion sheaf appearing on the left will
be zero, being a subsheaf of a torsion-free sheaf supported on a
codimension $1$ subvariety.
Then we tensor the short exact sequence with
$\pi_{\tC_1}^*\cO_{\{p_1,p_2\}}$ to obtain
$$ 0 \lfun T \lfun \cO_{\tC_1}(-p_i)|_{\{p_1,p_2\}} \boxtimes \pi_*
\cO_N \lfun \cO_{\{p_1,p_2\}} \boxtimes (\pi_p)_* \cO_{C_p}
\lfun \cO_{\{p_i\}} \boxtimes (\pi_p)_* \cO_{C_p} \lfun 0, $$
from which is clear that $T\simeq\cO_{\{p_i\}}\boxtimes(\pi_p)_*\cO_{C_p}$,
and therefore $\cF_1|_{\{p_i\}\times C}\simeq(\pi_p)_*\cO_{C_p}$.

An identical process defines a sheaf $\cF_2$ on $\tC_2\times C$ such that
\begin{eqnarray*}
\begin{CD}
\cE_{L,2}|_{\tC_2 \times C} @>>> \cA_2|_{\tC_2 \times C} \otimes \cF_2
@>>> 0 \\
@VVV @VVV \\
\cE_{L,1}|_{\tC_1 \times C} @>>> \cA_1|_{\tC_1 \times C} \otimes \cF_1
@>>>0
\end{CD},
\end{eqnarray*}
where we identify $\tC_2$ and $\tC_1$ via the isomorphism
$\varepsilon_2|_{\tC_2}$.
Since the cokernel of the vertical maps is again
$\cO_{\{\tp_1,\tp_2\}}\boxtimes(\pi_p)_*\cO_{C_p}$, the exact same proof as
above shows that $\cF_2|_{\{q\}\times C}\simeq\cO_C(q)$ if $q\in\tC_2$,
$q\neq\tp_1,\tp_2$, and $\cF_2|_{\{\tp_i\}\times C}\simeq(\pi_p)_*\cO_{C_p}$
for $i=1,2$.

We define $\cE_{L,3}$ to be the kernel of the map
$(\varepsilon_3,1)^*\cE_{L,2}\fun
(\varepsilon_3,1)^*(\cA_2|_{\tC_2\times C}\otimes\cF_2)$.

\section{Relation between $\cE_{L,3}$ and $\phi_{L,3}$}\label{rel3}

\begin{pro}\label{phiL3E3qxC}
For every $q\in\tC_2$ mapping to a smooth point $q\in C$,
$\phi_{L,3}|_{E_3|_q}$ is a morphism, and it maps $E_3|_q$ isomorphically to
$\bP(\eExt^1_C(L(-q),\cO_C(q)))$.
\end{pro}

\pf
We proved in Section \ref{third} that
$E_3|_q\simeq\bP(\Ext^1_C(L(-q),\cO_C(q)))$, and therefore we need to
show that under this identification, $\phi_{L,3}$ is the identity map.
The proposition follows from Proposition \ref{EL3E3qxC}.
\qed

\begin{pro}\label{EL3E3qxC}
The restriction of the torsion-free sheaf $\cE_{L,3}$ to $E_3|_q\times C$ is
the element of the vector space
$\eExt^1_{E_3|_q\times C}(\pi_C^*L(-q),\cO_{E_3|_q}(1)\boxtimes\cO_C(q))$
which corresponds to the identity under the identification of this extension
space with $\eHom(\eExt^1_C(L(-q),\cO_C(q)),\eExt^1_C(L(-q),\cO_C(q)))$.
\end{pro}

\pf
This proof is very similar to the proof of Proposition \ref{EL1E1xC}.
Let $H\su\bP_L$ be a linear hyperplane which contains $q$, does not
contain any node $p$ of $C$, and is transverse to the curve $C$ at $q$ (i.e.,
$H$ does not contain $T_qC$).
Then $H$ is isomorphic to its strict transform in $\bP_{L,2}$, that we
shall still denote by $H$.
It is clear that $\cE_{L,2}|_{H\times C}$ is $\cE_L|_{H\times C}$, and
therefore it is an extension
$ 0 \fun \pi_H^* \cO_H(1) \fun \cE_{L,2}|_{H \times C} \fun \pi_C^* L \fun 0 $
on $H\times C$.

Let $\sigma:\tH\fun H$ be the blow-up of $H$ at $q$, and let $E'\su\tH$
be the exceptional divisor.
Then there exists a commutative diagram on $\tH\times C$:
\begin{eqnarray*}
\begin{CD}
& & 0 & & 0 & & 0 \\
& & @VVV @VVV @VVV \\
0 @>>> \sigma^* \cO_H(1) \otimes \cO_{\tH}(-E') @>>>
\cE_{L,3}|_{\tH \times C} @>>> \cB_3|_{\tH \times C} @>>> 0 \\
& & @VVV @VVV @VVV \\
0 @>>> \sigma^* \cO_H(1) @>>> (\sigma,1)^* \cE_{L,2}|_{H \times C} @>>>
L @>>> 0 \\
& & @VVV @VVV @VVV \\
0 @>>> \cO_{E' \times C} @>>> \cO_{E'} \boxtimes \cO_C(q) @>>>
\cO_{E' \times \{q\}} @>>> 0 \\
& & @VVV @VVV @VVV \\
& & 0 & & 0 & & 0
\end{CD},
\end{eqnarray*}
where $\cB_3$ is defined exactly in the same way we defined $\cB_1$
and $\cB_2$, and the columns are exact (see Lemma \ref{codim}).

Let $\cE'_H$ be the push-forward of
$(\sigma,1)^*\cE_{L,2}|_{H\times C}$ via
$\sigma^*\cO_H(1)\inc\sigma^*\cO_H(1)\boxtimes\cO_C(q)$:
\begin{eqnarray}\label{cE'H}
\begin{CD}
0 @>>> \sigma^* \cO_H(1) @>>> (\sigma, 1)^* \cE_{L,2}|_{H \times C}
@>>> L @>>> 0 \\
& & @VVV @VVV || \\
0 @>>> \sigma^* \cO_H(1) \boxtimes \cO_C(q) @>>> \cE'_H @>>> L @>>> 0
\end{CD}.
\end{eqnarray}
Then the restriction of $\cE'_H$ to $E'\times C$ splits.
Indeed, $\sigma^*\cO_H(1)|_{E'}\simeq\cO_{E'}$, and via the identification
$\Ext^1_{E'\times C}(L,\cO_C(q))\simeq H^0(E',\cO_{E'})\otimes
\Ext^1_C(L,\cO_C(q))$ (see [Arc04]), we see that
$\cE'_H|_{E'\times C}$ splits as long as $\cE'_H|_{\{x\}\times C}$
splits for some $x\in E'$.
Resctricting the diagram (\ref{cE'H}) above to $\{x\}\times C$ for any
$x\in E'$, we see that $\cE'_H|_{\{x\}\times C}$ is the trivial
extension $\psi_q(\cE_L|_{\{q\}\times C})$.

Therefore, there exists a surjective map
$\cE'_H\fun\cO_{E'}\boxtimes\cO_C(q)$, and we can define $\cE'_{H,1}$
to be its kernel:
$ 0 \fun \cE'_{H,1} \fun \cE'_H \fun \cO_{E'} \boxtimes \cO_C(q) \fun 0. $
Then there exists the following commutative diagram on $\tH\times C$:
\begin{eqnarray*}
\begin{CD}
& & 0 & & 0 \\
& & @VVV @VVV \\
0 & \ \lfun\ & (\sigma^* \cO_H(1) \otimes \cO_{\tH}(-E')) \boxtimes
\cO_C(q) @>>> \cE'_{H,1} @>>> L @>>> 0 \\
& & @VVV @VVV || \\
0 & \lfun & \sigma^* \cO_H(1) \boxtimes \cO_C(q) @>>> \cE'_H @>>> L
@>>> 0\\
& & @VVV @VVV \\
& & \cO_{E'} \boxtimes \cO_C(q) & = & \cO_{E'} \boxtimes \cO_C(q) \\
& & @VVV @VVV \\
& & 0 & & 0
\end{CD}.
\end{eqnarray*}
Moreover, we have the following commutative diagram on $\tH\times C$
which relates $\cE'_H$ and $\cE'_{H,1}$ to $\cE_{L,2}|_{H\times C}$ and
$\cE_{L,3}|_{\tH\times C}$:
\begin{eqnarray*}
\begin{CD}
& & 0 & & 0 \\
& & @VVV @VVV \\
0 @>>> \cE_{L,3}|_{\tH \times C} @> i_1 >> (\sigma, 1)^*
\cE_{L,2}|_{H \times C} @>>> \cO_{E'} \boxtimes \cO_C(q) @>>> 0 \\
& & @VVV @VVV || \\
0 @>>> \cE'_{H,1} @> i'_1 >> \cE'_H @>>> \cO_{E'} \boxtimes \cO_C(q)
@>>> 0 \\
& & @VVV @VVV \\
& & \sigma^* \cO_H(1) \boxtimes \bC_q & = & \sigma^* \cO_H(1) \boxtimes
\bC_q \\
& & @VVV @VVV \\
& & 0 & & 0
\end{CD}.
\end{eqnarray*}
When we restrict the first two rows of this diagram to $E'\times C$,
and we look at the image of the restrictions of $i_1$ and $i'_1$ to
$E'\times C$, we obtain the following diagram, where the first row
shows that the restriction of $\cE_{L,3}$ to
$E_3|_q\times C\simeq E'\times C$ is an extension of $\pi_C^*L(-q)$ by
$\cO_{E_3|_q}(1)\boxtimes\cO_C(q)$:
\begin{eqnarray*}
\begin{CD}
0 @>>> \cO_{E'}(1) \boxtimes \cO_C(q) @>>> \cE_{L,3}|_{E' \times C}
@> i_1|_{E'\times C} >> L(-q) @>>> 0 \\
& & || & & @VVV @VVV \\
0 @>>> \cO_{E'}(1) \boxtimes \cO_C(q) @>>> \cE'_{H,1}|_{E' \times C}
@> i'_1|_{E'\times C} >> L @>>> 0 \\
\end{CD}.
\end{eqnarray*}
This shows that $\cE_{L,3}|_{E'\times C}$ is the pull-back of
$\cE'_{H,1}|_{E'\times C}$ via the pull-back of the inclusion
$L(-q)\inc L$ from $C$ to $E'\times C$.
Here is a summary of how to construct $\cE_{L,3}|_{E'\times C}$:
\begin{eqnarray*}
\begin{matrix}
\cE_{L,2}|_{H \times C} & \in & \Ext^1_{H \times C}(L, \cO_H(1)) \\
& & \downarrow \\
(\sigma,1)^* \cE_{L,2}|_{H \times C} & \in &
\Ext^1_{\tH \times C}(L, \sigma^* \cO_H(1)) \\
& & \downarrow \\
\cE'_H & \in & \Ext^1_{\tH \times C}(L, \sigma^* \cO_H(1) \boxtimes
\cO_C(q)) \\
& & \uparrow \\
\cE'_{H,1} & \in & \Ext^1_{\tH \times C}(L, (\sigma^* \cO_H(1) \otimes
\cO_{\tH}(- E')) \boxtimes \cO_C(q))\\
& & \downarrow \\
\cE'_{H,1}|_{E' \times C} & \in & \Ext^1_{E' \times C}(L, \cO_{E'}(1)
\boxtimes \cO_C(q)) \\
& & \downarrow \\
\cE_{L,3}|_{E' \times C} & \in &
\Ext^1_{E' \times C}(L(-q), \cO_{E'}(1) \boxtimes \cO_C(q))
\end{matrix}
\end{eqnarray*}
Using the isomorphisms
$ \Ext^1_{Y \times C} (L, F \boxtimes G) \simeq H^0(Y, F) \otimes
\Ext^1_C(L, G) $
(see [Arc04]), we can understand what extension
$\cE_{L,3}|_{E'\times C}$ is by tracking the corresponding elements in
these spaces.
Let $v_0,\dots,v_n$ be a basis of $\Ext^1_C(L,\cO_C)$ with
$\Span\{v_1,\dots,v_n\}=\langle H\rangle$, and
$\Span\{v_0,v_1\}=\langle T_qC\rangle$.
Let $v_0^*,\dots,v_n^*$ be the corresponding dual basis in
$\Ext^1_C(L,\cO_C)^*$.
Then $v_1^*,\dots,v_n^*$ is a basis of
$\langle H\rangle^*\simeq H^0(H,\cO_H(1))$, and $\cE_{L,2}|_{H\times C}$
corresponds to the element
$ \sum_{i=1}^n v_i^* \otimes v_i \in H^0(H, \cO_H(1)) \otimes
\Ext^1_C(L, \cO_C). $
Let $\psi_{-q}\colon\Ext^1_C(L,\cO_C(q))\fun\Ext^1_C(L(-q),\cO_C(q))$ be the
natural linear homomorphism.
Since $\psi_q(v_1)=0$ and $\ker(\psi_{-q}\circ\psi_q)=\Span\{v_0,v_1\}$, we
can calculate that $\cE_{L,3}|_{E'\times C}$ corresponds to the element
$\sum_{i=2}^n w_i^*\otimes w_i\in H^0(E',\cO_{E'}(1))\otimes
\Ext^1_C(L(-q),\cO_C(q))$, where, for $2\leq i\leq n$,
$w_i=\psi_{-q}(\psi_q(v_i))$.
Therefore, $\cE_{L,3}$ corresponds to the identity in the vector space
$\Hom(\Ext^1_C(L(-q),\cO_C(q)),\Ext^1_C(L(-q),\cO_C(q)))$, as claimed.
\qed

We now prove that, for every node $p$ of $C$, and for every $i\in\{1,2\}$,
$\phi_{L,3}|_{\bP(T_{\tp_i}(E_2|_{p_i}))}$ factors through the canonical
isomorphism
$\bP(T_{\tp_i}(E_2|_{p_i}))\overset{\simeq}{\fun}\bP(\Im\psi_{L_p})$
described in Section \ref{third}.

\begin{pro}\label{EL3Impsi}
For every node $p$ in $C$, and every $i\in\{1,2\}$, the extension
$$ \cE_{L,3}|_{\bP(\cN_{\{\tp_i\}/E_2|_{p_i}}) \times C} \in
\eExt^1_{\bP(\cN_{\{\tp_i\}/E_2|_{p_i}}) \times C}(L \otimes
((\pi_p)_*\cO_{C_p})^*, \cO_{\bP(\cN_{\{\tp_i\}/E_2|_{p_i}})}(1) \boxtimes
(\pi_p)_*\cO_{C_p}) $$
corresponds to the inclusion in
$\eHom(\emph{Im}\psi_{L_p},
\eExt^1_C(L\otimes((\pi_p)_*\cO_{C_p})^*,(\pi_p)_*\cO_{C_p}))$
under the canonical identification
$\cN_{\{\tp_i\}/E_2|_{p_i}}\simeq\emph{Im}\,\psi_{L_p}$.
\end{pro}

\pf
Let $i\in\{1,2\}$.
From our description of $\cE_{L,2}$ and $\phi_{L,2}$, it is clear that
$\cE_{L,2}|_{E_2|_{p_i}\times C}$ induces the linear map given by projection
from $\tp_i$, i.e., it corresponds to a linear homomorphism
$\cN_{L_p/\bP_{L,1}}|_{p_i}\fun
\Ext^1_C(L\otimes((\pi_p)_*\cO_{C_p})^*,(\pi_p)_*\cO_{C_p})$ of kernel
$\langle\tp_i\rangle$ and image $\Im\psi_{L_p}$.
Therefore, we can find a basis $w_1,\dots,w_n$ of
$\cN_{L_p/\bP_{L,1}}|_{p_i}$ and a basis $v_0,\dots,v_n$ of
$\Ext^1_C(L\otimes((\pi_p)_*\cO_{C_p})^*,(\pi_p)_*\cO_{C_p})$ such that
$\langle\tp_i\rangle=\Span\{w_1\}$, $\Im\psi_{L_p}=\Span\{v_2,\dots,v_n\}$,
and $\Span\{w_1,w_j\}$ maps to $\Span\{v_j\}$ for every $2\leq j\leq n$
under the homomorphism corresponding to
$\cE_{L,2}|_{E_2|_{p_i}\times C}$.
In particular, this sheaf corresponds to $\sum_{j=2}^nw_j^*\otimes v_j$
in $H^0(E_2|_{p_i},\cO_{E_2|_{p_i}}(1))\otimes
\Ext^1_C(L\otimes((\pi_p)_*\cO_{C_p})^*,(\pi_p)_*\cO_{C_p})$.

To simplify the notation, let us denote $E_2|_{p_i}$ by $X$ and its
blow-up at $\tp_i$ by $\sigma:\tX\fun X$, with $E'$ the exceptional
divisor.
Then there exists a short exact sequence
$$ 0 \lfun \cE_{L,3}|_{\tX \times C} \lfun (\sigma,1)^*
\cE_{L,2}|_{X \times C} \lfun \cO_{E'} \boxtimes (\pi_p)_*\cO_{C_p} \lfun 0,$$
obtained by restricting the short exact sequence defining $\cE_{L,3}$
to $\tX\times C$.
It stays exact because of Lemma \ref{codim}.

There exists the following commutative diagram
\begin{eqnarray*}
\begin{CD}
& & 0 & & 0 \\
& & @VVV @VVV \\
0 & \ \lfun \ & \cA \boxtimes (\pi_p)_* \cO_{C_p} @>>>
\cE_{L,3}|_{\tX \times C} @>>> L \otimes ((\pi_p)_* \cO_{C_p})^* @>>> 0 \\
& & @VVV @VVV || \\
0 & \lfun & \sigma^* \cO_X(1) \boxtimes (\pi_p)_* \cO_{C_p} @>>> (\sigma, 1)^*
\cE_{L,2}|_{X \times C} @>>> L \otimes ((\pi_p)_* \cO_{C_p})^* @>>> 0 \\
& & @VVV @VVV \\
& & \cO_{E'} \boxtimes (\pi_p)_*\cO_{C_p} & = & \cO_{E'} \boxtimes
(\pi_p)_*\cO_{C_p}\\
& & @VVV @VVV \\
& & 0 & & 0
\end{CD},
\end{eqnarray*}
where $\cA=\sigma^*\cO_X(1)\otimes\cO_\tX(-E')$.
If we restrict the first row to $E'\times C$, we obtain a short exact
sequence
$$ 0 \lfun \cO_{E'}(1) \boxtimes (\pi_p)_* \cO_{C_p} \lfun
\cE_{L,3}|_{E' \times C} \lfun L \otimes ((\pi_p)_* \cO_{C_p})^* \lfun 0. $$
Remember that $E'$ is $\bP(\cN_{\{\tp_i\}/E_2|_{p_i}})$.
The following diagram, where, to simplify the notation, we denoted
$\Ext^1_C(L\otimes((\pi_p)_*\cO_{C_p})^*,(\pi_p)_*\cO_{C_p})$ by $V$,
illustrates the steps we took in finding $\cE_{L,3}|_{E'\times C}$:
\begin{eqnarray*}
\begin{matrix}
\cE_{L,2}|_{X \times C} & \longleftrightarrow & \ds\sum_{j=2}^n w_j^*
\otimes v_j & \in & H^0(X, \cO_X(1)) \otimes V \\
& & & & \downarrow \\
(\sigma,1)^* \cE_{L,2}|_{X \times C} & \longleftrightarrow &
\ds\sum_{j=2}^n w_j^* \otimes v_j & \in & H^0(\tX, \sigma^* \cO_X(1))
\otimes V \\
& & & & \uparrow \\
\cE_{L,3}|_{\tX \times C} & \longleftrightarrow & \ds\sum_{j=2}^n w_j^*
\otimes v_j & \in & H^0(\tX, \cA) \otimes V \\
& & & & \downarrow \\
\cE_{L,3}|_{E' \times C} & \longleftrightarrow & \ds\sum_{j=2}^n v_j^*
\otimes v_j & \in & H^0(E', \cO_{E'}(1)) \otimes V
\end{matrix}.
\end{eqnarray*}
Therefore, $\cE_{L,3}|_{E'\times C}$ corresponds to the inclusion
$\Im\psi_{L_p}\inc
\Ext^1_C(L\otimes((\pi_p)_*\cO_{C_p})^*,(\pi_p)_*\cO_{C_p})$, as claimed.
\qed

Let $p$ be a node of $C$, let $l\in L_p$, $l\neq p_1,p_2$, and let
$Y_l:=s_{E_l}(L_p)$, where $s_{E_l}$ is the section of $E_{2,p}\fun L_p$
defined in Section \ref{third}.
Remember that we denoted by $y_{l,i}$ the only point of intersection of the
strict transform $\tY_l$ of $Y_l$ with $E_3|_{\tp_i}$ ($i=1,2$).

\begin{pro}\label{E3l}
The restriction of $\cE_{L,3}$ to $\tY_l\times C$ is a non-zero element in the
vector space $\eExt^1_{\tY_l\times C}(L\otimes((\pi_p)_*\cO_{C_p})^*,
(\pi_p)_*\cO_{C_p})\simeq H^0(\tY_l,\cO_{\tY_l})\otimes
\eExt^1_C(L\otimes((\pi_p)_*\cO_{C_p})^*,(\pi_p)_*\cO_{C_p})$.
In particular $\cE_{L,3}|_{\{y_{l,i}\}\times C}\simeq E_l$ for $i=1,2$.
\end{pro}

For the proof, we need the following result.

\begin{lem}
The restriction of $\cO_{\bP_{L,2}}(-E_2)$ to $L_{2,p}\simeq\bP^1$ is
$\cO_{L_{2,p}}(1)$.
\end{lem}

\pf
Remember that $L_{2,p}=\widetilde{T_pC}\cap E_2$.
It is isomorphic to $L_p=\widetilde{T_pC}\cap E_1$ via $\varepsilon_2$,
and it is therefore the exceptional divisor of the blow-up of $T_pC$ at $p$.
Therefore,
$$ \cO_{\bP_{L,2}}(-E_2)|_{L_{2,p}} =
(\cO_{\bP_{L,2}}(-E_2)|_{\widetilde{T_pC}})|_{L_{2,p}} =
\cO_{\widetilde{T_pC}}(- L_{2,p})|_{L_{2,p}} = \cO_{L_{2,p}}(1). $$
\qed

{\em Proof $($of Proposition $\ref{E3l})$.}
We saw in Lemma \ref{cEL2E2} that there exists a short exact sequence
$ 0 \fun (\varepsilon_2^* \cO_{L_p}(1) \otimes \cO_{E_2}(-E_2))
\boxtimes (\pi_p)_* \cO_{C_p} \fun \cE_{L,2}|_{E_2 \times C} \fun L \otimes
((\pi_p)_* \cO_{C_p})^* \fun 0 $
on $E_2\times C$.
If we restrict it to $Y_l\times C$, we obtain the short exact
sequence
$$ 0 \lfun \cO_{Y_l}(2) \boxtimes (\pi_p)_* \cO_{C_p} \lfun
\cE_{L,2}|_{Y_l \times C} \lfun L \otimes ((\pi_p)_* \cO_{C_p})^* \lfun 0, $$
because $\cO_{E_2}(-E_2)|_{Y_l}=\cO_{Y_l}(1)$ since $Y_l$ and $L_2$ are
in the same linear system.
Therefore, there exists the following commutative diagram on
$\tY_l\times C\simeq Y_l\times C$:
\begin{eqnarray*}
\begin{CD}
& & 0 & & 0 \\
& & @VVV @VVV \\
0 & \ \lfun \ & \cO_{\tY_l \times C} \boxtimes (\pi_p)_* \cO_{C_p} @>>>
\cE_{L,3}|_{\tY_l \times C} @>>> L \otimes ((\pi_p)_* \cO_{C_p})^* @>>> 0 \\
& & @VVV @VVV || \\
0 & \lfun & \cO_{Y_l \times C}(2) \boxtimes (\pi_p)_* \cO_{C_p} @>>>
\cE_{L,2}|_{Y_l \times C} @>>> L \otimes ((\pi_p)_* \cO_{C_p})^* @>>> 0 \\
& & @VVV @VVV \\
& & \cO_{\{\tp_1,\tp_2\}} \boxtimes (\pi_p)_* \cO_{C_p} & = &
\cO_{\{\tp_1,\tp_2\}} \boxtimes (\pi_p)_* \cO_{C_p} \\
& & @VVV @VVV \\
& & 0 & & 0
\end{CD},
\end{eqnarray*}
and $\cE_{L,3}|_{\tY_l \times C}$ is an element of
$\Ext^1_{\tY_l\times C}(L\otimes((\pi_p)_*\cO_{C_p}),(\pi_p)_*\cO_{C_p})$.
Since this is isomorphic to $H^0(\tY_l,\cO_{\tY_l})\otimes
\Ext^1_C(L\otimes((\pi_p)_*\cO_{C_p}),(\pi_p)_*\cO_{C_p})$ (see [Arc04]),
$\cE_{L,3}|_{\tY_l \times C}$ is non-zero because we know that
$\cE_{L,3}|_{\{y\}\times C}$ does not split for every $y\in\tY_l$,
$y\neq y_{l,1},y_{l,2}$.
\qed

\section{Fibers of $\phi_{L,3}$}\label{fibers}

We prove in this section that the fibers of $\phi_{L,3}$ are connected.
Let us start with characterizing the image.

\begin{lem}\label{image}
An element $E\in\oSU$ is in the image of $\phi_{L,3}$ if and only if
$H^0(E)\neq0$.
\end{lem}

\pf
By our description of $\phi_{L,3}$ it is clear that if $E$ is in its image,
then there exists a non-zero map $\cO_C\fun E$, and therefore $E$ has a
non-zero section.
Conversely, if $H^0(E)\neq0$, there exists a non-zero section
$\cO_C\overset{s}{\fun}E$.
If $E$ is stable, then $s$ can vanish at at most one point, and therefore $E$
fits in at least one of the following exact sequences:
\begin{itemize}
\item
$ 0 \fun \cO_C \fun E \fun L \fun 0, $
\item
$ 0 \fun \cO_C(q) \fun E \fun L(-q) \fun 0 $ for some smooth point $q\in C$,
\item
$ 0 \fun (\pi_p)_* \cO_{C_p} \fun E \fun L \otimes ((\pi_p)_* \cO_{C_p})^*
\fun 0 $ for some node $p\in C$,
\end{itemize}
and it is therefore in the image of $\phi_{L,3}$.
This proves the lemma, since $\phi_{L,3}$ is proper and the locus of stable
bundles is dense in $\{E\in\oSU\mid H^0(E)\neq0\}$.
\qed

As in the smooth case, we have the following result.

\begin{pro}\label{h0E}
For every stable $E\in\oSU$ there exists a morphism
$$ \psi_E \colon \bP(H^0(E)) \lfun \bP_L $$
such that, for every $x\in\bP_L\setminus C$,
$$ \phi_L(x) = E \Longleftrightarrow x \in \emph{Im}(\psi_E). $$
\end{pro}

\pf
For every $s\in H^0(E)$, define
$ \psi_E([s]) := \bP(\ker(\Ext^1_C(L, \cO_C) \fun \Ext^1_C(L, E))). $
To prove that $\psi_E$ is a morphism, it suffices to show that, for every
$\in H^0(E)$, the kernel in the definition of $\psi_E([s])$ is
one-dimensional.
Note that, since $E$ is stable of degree $\leq4$, every torsion-free subsheaf
of $E$ of rank $1$ has degree $\leq1$, and therefore $s$ can vanish at at most
one point.

Case I: $s$ is no-where vanishing.
There exists a short exact sequence
$ 0 \fun \cO_C \fun E \fun L \fun 0, $
and applying to it the functor $\Hom_C(L,-)$ we obtain
$$ 0 \lfun \Hom_C(L, L) \lfun \Ext^1_C(L, \cO_C) \lfun \Ext^1_C(L, E) \lfun
\cdots, $$
where the sequence starts with $\Hom_C(L,E)$, which is zero because $E$ is
stable.
This proves that the kernel of $\Ext^1_C(L,\cO_C)\fun\Ext^1_C(L,E)$ is
isomorphic to $\Hom_C(L,L)$, and it is therefore one-dimensional.
Moreover, the image of the identity element of $\Hom_C(L,L)$ in
$\Ext^1_C(L,\cO_C)$ is the extension associated to $\cO_C\overset{s}{\fun}E$,
and therefore $\phi_L(\psi_E([s]))=E$.

Case II: $s$ vanishes at exactly one point.
Let $F$ be $\cO_C(q)$ if $s$ vanishes at a smooth point $q\in C$ or let $F$ be
$(\pi_p)_*\cO_{C_p}$ if $s$ vanishes at a node $p$ of $C$.
There exists a short exact sequence
\begin{eqnarray}\label{H^0E}
0 \lfun F \lfun E \lfun L \otimes F^* \lfun 0,
\end{eqnarray}
and $\Ext^1_C(L,\cO_C)\fun\Ext^1_C(L,E)$ factors through 
$\Ext^1_C(L,F)\fun\Ext^1_C(L,E)$.
Since the kernel of $\Ext^1_C(L,\cO_C)\fun\Ext^1_C(L,F)$ is one-dimensional,
to conclude the proof it suffices to show that
$\Ext^1_C(L,F)\fun\Ext^1_C(L,E)$ is injective.
Applying the functor $\Hom_C(L,-)$ to (\ref{H^0E}), we see that this is the
case, because $\Hom_C(L,L\otimes F^*)=0$.
\qed

In what follows, we also need the following result, which is similar to the
proposition above.

\begin{lem}\label{h0Ept}
Let $E\in\oSU$ be stable.

$(a)$ If every section $s\in H^0(E)$ vanishes at a smooth point $q\in C$, then
there exists a morphism
$$ \psi_{E,q} \colon \bP(H^0(E)) \lfun E_3|_q \simeq
\bP(\eExt^1_C(L(-q), \cO_C(q))) $$
such that
$$ \phi_{L,3}^{-1}(E) = \eIm \phi_{E,q}. $$

$(b)$ If every section $s\in H^0(E)$ vanishes at a node $p\in C$, then there
exists a morphism
$$ \psi_{E,p} \colon \bP(H^0(E)) \lfun \bP(H'_p) \su
\bP(\eExt^1_C(L \otimes ((\pi_p)_* \cO_{C_p})^*, (\pi_p)_* \cO_{C_p})) $$
such that $x$ is in the image of $\psi_{E,p}$ if and only if $x$ maps to $E$
under the natural forgetful map $\bP(H'_p)\fun\oSU$, which is a morphism by
Lemma $\ref{extss}$.
\end{lem}

\pf
Define
$\psi_{E,q}([s]):=\bP(\ker(\Ext^1_C(L(-q),\cO_C(q))\fun\Ext^1_C(L(-q),E)))$ in
part (a), and $\psi_{E,p}([s]):=\bP(\ker(\Ext^1_C(L\otimes
((\pi_p)_*\cO_{C_p})^*,(\pi_p)_*\cO_{C_p})\fun\Ext^1_C(L\otimes
((\pi_p)_*\cO_{C_p})^*,E)))$ in part (b).
If we let $F=\cO_C(q)$ for part (a) and $F=(\pi_p)_*\cO_{C_p}$ for part (b),
the proof is the same as the proof of Case I in Proposition \ref{h0E} where we
now use the unique extension $0\fun F\fun E\fun L\otimes F^*\fun0$ associated
to $s$.
\qed

To simplify the notation, for the rest of this section, we shall say that an
element $E\in\oSU$ is of type Q if there exists a non-zero map $\cO_C(q)\fun
E$ for some smooth point $q\in C$, and is of type P if there exists a non-zero
map $(\pi_p)_*\cO_{C_p}\fun E$ for some node $p\in C$.

\begin{pro}\label{fiberspro}
The fibers of $\phi_{L,3}$ are connected.
\end{pro}

We shall divide the proof of this proposition into several lemmas analyzing
various different cases.

\begin{lem}\label{fiber1}
If $g>\deg L$, and $E\in\eIm\phi_{L,3}$ is not of type P, then
$\phi_{L,3}^{-1}(E)$ is just a point.
\end{lem}

\pf
Since $E$ is not of type P, there are two possibilities:

Case I:
There exists an $x\in\bP_L\setminus C$ such that $\phi_L(x)=E$.
Then there exists a short exact sequence $0\fun\cO_C\fun E\fun L\fun0$, and
$h^0(E)=h^0(\cO_C)=1$, since $H^0(L)=0$ because $g>\deg L$ and $L$ is generic.
This proves that there is only one way to write $E$ as an extension of $L$ by
$\cO_C$, and $\cO_C$ must be a maximal subbundle of $E$ (i.e., $E$ is not in
the image of the natural morphisms $\bP(\Ext^1_C(L(-q),\cO_C(q)))\fun\oSU$ and
$\bP(H'_p)\lfun\oSU$).
Therefore $\phi_{L,3}^{-1}(E)=\{pt\}$.

Case II:
There exists a smooth point $q\in C$ and an $x\in E_3|_q$ such that
$\phi_{L,3}(x)=E$.
Then there exists a short exact sequence $0\fun\cO_C(q)\fun E\fun L(-q)\fun0$,
and again $h^0(E)=h^0(\cO_C(q))=1$ and $\phi_{L,3}^{-1}(E)=\{pt\}$.
\qed

\begin{lem}\label{fiber2}
If $g>\deg L$, and $E\in\eIm\phi_{L,3}$ is a non-locally-free sheaf of type P,
then $\phi_{L,3}^{-1}(E)$ is the union of $($the strict transforms of\,$)$ a
plane and two lines intersecting the plane.
\end{lem}

\pf
If $E$ is of type P, then $E$ is in the image of a point of $\tE_1$, $\tE_2$,
or $E_3|_{\tp_i}$ for some node $p$ of $C$ and $i\in\{1,2\}$.
There exists a short exact sequence
$ 0 \fun (\pi_p)_* \cO_{C_p} \fun E \fun L \otimes ((\pi_p)_* \cO_{C_p})^*
\fun 0. $
Since $L$ is generic, $H^0(L)=0$, and this implies that
$H^0(L\otimes((\pi_p)_*\cO_{C_p})^*)=0$.
Then $h^0(E)=h^0((\pi_p)_*\cO_{C_p})=1$, and there is only one way to write
$E$ as such an extension (and $E$ cannot be written as an extension of $L$ by
$\cO_C$ or $L(-q)$ by $\cO_C(q)$ for some smooth $q\in C$ or
$L\otimes((\pi_{p'})_*\cO_{C_{p'}})^*)$ by $(\pi_{p'})_*\cO_{C_{p'}}$ for some
other node $p'$).

Since $E$ is not locally-free, it is in the image of $\psi_{L_p}$, and by
Theorem \ref{1} $\phi_{L,1}^{-1}(E)$ is a plane in the projective space
$E_1|_p$ containing the line $L_p$ (except for the line itself, where
$\phi_{L,1}$ is not defined).
Then, by Theorem \ref{2}, $\phi_{L,2}^{-1}(E)$ is the union of the strict
transform $\wt{\ov{\phi_{L,1}^{-1}(E)}}$ of the plane
$\ov{\phi_{L,1}^{-1}(E)}$ and two lines which are contained in $E_2|_{p_1}$
and $E_2|_{p_1}$, respectively (except for the points $\tp_1$ and $\tp_2$,
which are on the lines, where $\phi_{L,2}$ is not defined).
The lines intersect the plane at the points
$\wt{\ov{\phi_{L,1}^{-1}(E)}}\cap E_2|_{p_1}$ and
$\wt{\ov{\phi_{L,1}^{-1}(E)}}\cap E_2|_{p_1}$, respectively.
The last blow-up just adds the two missing points, and $\phi_{L,3}^{-1}(E)$ is
the union of a plane and two lines which intersect it, as claimed.
\qed

\begin{lem}\label{fiber3}
If $g>\deg L$, and $E\in\eIm\phi_{L,3}$ is a locally-free sheaf of type P,
then $\phi_{L,3}^{-1}(E)$ is $($the strict transform of\,$)$ a line.
\end{lem}

\pf
The proof of this lemma follows exactly the proof of the previous lemma up to
the description of $\phi_{L,1}^{-1}(E)$.
In our case now, since $E$ is locally-free, $\phi_{L,1}^{-1}(E)$ is empty.
Then, by Theorem \ref{2}, $\phi_{L,2}^{-1}(E)$ is a section of
$E_{2,p}\fun L_p$ which passes through the points $\tp_1$ and $\tp_2$, where
$\phi_{L,2}$ is not defined.
The last blow-up just adds the two missing points, and therefore
$\phi_{L,3}^{-1}(E)$ is isomorphic to the line $L_p$, as claimed.
\qed

The previous lemmas prove that the fibers of $\phi_{L,3}$ are connected if
$g>\deg L$.
To prove that the fibers are always connected, we need to still study the
cases of $g=2$, $g=3$ and the case $g=\deg L=4$.
We shall now prove the case $g=2$ and $\deg L=3$ (since it is the case for
which we have an application), and leave the other cases as an exercise for
the reader.

\begin{lem}
If $g=2$ and $\deg L=3$, then the fibers of $\phi_{L,3}$ are connected.
\end{lem}

\pf
For every $E\in\oSU$, $E$ is stable, and $H^0(E)\neq0$ because $\chi(E)=1$.
Therefore, $\phi_{L,3}$ is surjective by Lemma \ref{image}, and it is a
birational morphism.
If $h^0(E)=1$, then the fiber $\phi_{L,3}^{-1}(E)$ is the same as the fibers
described in Lemmas \ref{fiber1}, \ref{fiber2}, and \ref{fiber3} above.
Suppose that $h^0(E)\geq2$.
Since $E$ is stable, every section $s\in H^0(E)$ satisfies $\deg(Z(s))\leq1$.
Moreover, we have the morphism $\psi_E\colon\bP(H^0(E))\fun\bP_L$ described in
Proposition \ref{h0E}, and we saw that $\phi_L^{-1}(E)=\Im\psi_E\setminus C$.

Case I: $\Im\psi_E\not\su C$.
If $\Im\psi_E$ does not intersect $C$, then
$\phi_{L,3}^{-1}(E)=\wt{\Im\psi_E}$ is connected.
If $\Im\psi_E$ intersects $C$ at a smooth point $q$, then there exists a
unique section $s\in H^0(E)$ such that $Z(s)=\{q\}$, and $E$ can be written as
an extension of $L(-q)$ by $\cO_C(q)$.
By continuity, this extension must be the point in
$E_3|_q\simeq\bP(\Ext^1_C(\cO_C(q),L(-q))$ in the strict transform of the
closure of $\Im\psi_E$.
Similarly, if $\Im\psi_E$ intersects $C$ at a node $p$, then $E$ can be
written as an extension of $L\otimes((\pi_p)_*\cO_{C_p})^*$ by
$(\pi_p)_*\cO_{C_p}$ in a unique way.
Since $\Im\psi_E\not\su C$, $E$ is locally-free.
Therefore, $E$ is not in the image of $E_1$, and its preimage in $\wt E_2$ is
a line, which, by continuity, must intersect the strict transform of the
closure of $\Im\psi_E$.

Case II: $\Im\psi_E\su C$.
In this case, $\Im\psi_E$ must be just a point.
Indeed, $\psi_E(s)=x\in C$ if and only if $Z(s)=\{x\}$, and if $\Im\psi_E=C$,
then there would exist two distinct sections in $H^0(E)$ mapping to each node
of $C$.
If $\Im\psi_E$ is a smooth point $q$ of $C$, then every section of $E$
vanishes at $q$, and by Lemma \ref{h0Ept} $\phi_{L,3}^{-1}(E)=\Im\psi_{E,q}$
is connected.
If $\Im\psi_E$ is a node $p$ of $C$, then every section of $E$ vanishes at
$p$, and by Lemma \ref{h0Ept} there exists a morphism
$$ \psi_{E,p} \colon \bP(H^0(E)) \lfun \bP(H'_p) \su
\bP(\Ext^1_C(L\otimes((\pi_p)_*\cO_{C_p})^*, (\pi_p)_*\cO_{C_p})) $$
such that $x$ is in the image of $\psi_{E,p}$ if and only if $x$ maps to $E$
under the natural forgetful morphism $\bP(H'_p)\fun\oSU$.
For each point $x$ of $\bP(H'_p)$, the space of all points in
$\wt{E_{2,p}}\cup E_3|_{\tp_1}\cup E_3|_{\tp_2}$ which map to $x$ under the
maps to $\bP(H'_p)$ of Theorems \ref{2} and \ref{3} is connected by Lemmas
\ref{fiber2} and \ref{fiber3}.
Since the image of $\psi_{E,p}$ in $\bP(H'_p)$ is connected, this proves that
$\phi_{L,3}^{-1}(E)$ is also connected.
\qed

\begin{rem}
The proof shows that $\phi_{L,3}^{-1}(E)$ is connected for every stable $E$ in
the image of $\phi_{L,3}$.
\end{rem}

\section{The case $\deg L>4$}

If $\deg L>4$, the rational map $\bP_{L,3}\fun\oSU$ is not a morphism, but we
can still prove the following.

\begin{pro}\label{UVproper}
Let $\deg L=2g-1$, and let
$$ V = \{ E \in \oSU \mid H^0(E) \simeq \bC \cdot s \textrm{ and }
Z(s) = \emptyset \textrm{ or } \{q\}, \textrm{ with } q \textrm{ a smooth
point of } C \}. $$
Then the codimension of $\oSU\setminus V$ in $\oSU$ is $\geq2$, and there
exists an open subset $U$ of $\bP_{L,3}$ such that
$\phi_{L,3}|_U\colon U\fun V$ is an isomorphism.
\end{pro}

\pf
Let $E\in V$.
If $Z(s)=\emptyset$ [resp.\ $Z(s)=\{q\}$], then $E$ is an extension of $L$ by
$\cO_C$ [resp.\ of $L(-q)$ by $\cO_C(q)$], and by Lemma \ref{h0E} [resp.\
\ref{h0Ept}] there exists a unique point $x$ of $\bP_L$ [resp.\ $\bP_{L,3}$]
such that $\phi_L(x)=E$ [resp.\ $\phi_{L,3}(x)=E$].
No other point of $\bP_{L,3}$ can map to $E$, and therefore $E$ has a unique
preimage under $\phi_{L,3}$.

To prove the claim about the codimension of $\oSU\setminus V$ in $\oSU$, let
us first study $U$.
If we identify $\bP_L\setminus C$ with its isomorphic image in $\bP_{L,3}$,
then
$$ U \cap (\bP_L \setminus C) = \{ E \in \bP_L \mid E \textrm{ is semi-stable
and } h^0(E) = 1 \}. $$
As in the smooth case, let
$$ \Gamma_L := \{ E \in \bP_L \mid h^0(E) > 1 \}. $$
It is a hypersurface of degree $g$ in $\bP_L$ (see [Ber92]).
Then
$$ \bP_L = (U \cap (\bP_L \setminus C)) \cup \Gamma_L \cup B, $$
where $B$ is the base locus of $\phi_L$, which has codimension $\geq2$ in
$\bP_L$.
Moreover, $U$ does not intersect the exceptional divisors $E_1$ and $E_2$, and
$U\cap E_3$ is a dense open subset of $E_3$.
Therefore, the complement of $U$ in $\bP_{L,3}$ is the union of $\tG_L$,
$E_1$, $E_2$, and a locus of codimension $\ge2$ in $\bP_{L,3}$.
Since the map $\phi_{L,3}$ restricted to $\tG_L\cup E_1\cup E_2$ has positive
dimensional fibers, the image of $\tG_L\cup E_1\cup E_2$ in $\oSU$ has
codimension $2$, and therefore
$ \codim(\oSU \setminus V, \oSU) \geq 2. $
\qed

\section{Applications}\label{applications}

Before we give direct applications of our construction, let us point out how
the rational map $\phi_L$ can be used to describe $\oSU$ on an irreducible
nodal curve of genus $1$.
In this case, the normalization $N$ of $C$ is isomorphic to $\bP^1$, and $C$
has only one node.

\begin{pro}
Let $C$ be an irreducible projective curve of arithmetic genus $1$ with one
node $p$ as singularity.

$(1)$ Let $L$ be any line bundle of degree $1$.
Then
$$ \phi_L \colon \bP(\eExt^1_C(L, \cO_C)) \lfun \SU_C(2,L) \su \oSU $$
is an isomorphism, and therefore $\oSU=\SU_C(2,L)\simeq\bP^0$ as in the smooth
case $($see $[\emph{Tu}93])$.

$(2)$ Let $L$ be any line bundle of degree $2$.
Then
$$ \phi_L \colon \bP(\eExt^1_C(L, \cO_C)) \lfun \SU_C(2,L) \su \oSU $$
is an isomorphism, and therefore $\oSU=\SU_C(2,L)\simeq\bP^1$ as in the smooth
case $($see $[\emph{Tu}93])$.
\end{pro}

\pf
The map $\phi_L$ is a morphism by Proposition \ref{UL}.
Note that every $E\in\oSU$ has at least one section because
$h^0(E)\geq\chi(E)=\deg L\geq1$.

If $\deg L=1$, since every element of $\oSU$ is stable, the sections of $E$
cannot vanish at any point, therefore $E$ is an extension of $L$ by $\cO_C$,
and $\phi_L$ is surjective.

If $\deg L=2$, then $h^0(L)=2$.
We have that $\Ext^1_C(L,\cO_C)\simeq H^0(L)^*$, and therefore
$\bP_L\simeq|L|^*$, which, being a $\bP^1$, is canonically equal to $|L|$.
The morphism $\phi_L$ is defined by
$\phi_L(q_1+q_2)=\cO_C(q_1)\oplus\cO_C(q_2)$ for
$q_1+q_2\in|L|\simeq\bP(\Ext^1_C(L,\cO_C))$.
It is clearly injective.
To prove that $\phi_L$ is surjective, note that, if $E\in\oSU$ is not stable,
then it must S-equivalent to $\cO_C(q_1)\oplus\cO_C(q_2)$ with
$q_1+q_2\in|L|$, and if it is stable, then the sections of $E$ cannot vanish
at any point, and $E$ is in the image of $\phi_L$.
\qed

From the direct description of $\phi_L$ in the proof, we deduce the following
fact, which is also true in the smooth case (see [Tu93]).

\begin{cor}
If $C$ is an irreducible projective curve of arithmetic genus $1$ with one
node $p$ as singularity, and $L$ is a line bundle of even degree, then every
vector bundle $E\in\SU_C(2,L)$ is semi-stable but not stable.
\end{cor}

For an irreducible nodal curve of genus $2$, we have the following application
as a corollary of Proposition \ref{fiberspro}.
This fact is already known (see [BhoNew90]).

\begin{cor}
If $g=2$ and $\deg L$ is odd, the normalization morphism $\oSU^\nu\fun\oSU$ is
one-to-one.
\end{cor}

\pf
This follows from the fact that $\phi_{L,3}$ is a birational morphism with
connected fibers for $g=2$ and $\deg L=3$.
\qed

\begin{rem}
If $g=2$ and $\deg L=4$, then $\phi_L\colon\bP_L\fun\oSU$ is a rational map
defined by sections of $|\cI_C(2)|$.
It should be possible to prove that $h^0(\cI_C(2))=4$, and obtain as a
corollary the known fact that, if $\deg L$ is even, $\oSU\simeq\bP^3$ for an
irreducible nodal curve of genus $2$ (see [Bho98]).
We can prove that $h^0(\cI_C(2))$ is indeed $4$ for a generic such curve with
one node and for a generic such curve with two nodes.
\end{rem}

For curves of genus $\geq2$, as a corollary of Proposition \ref{UVproper}, we
can prove the following results.

\begin{cor}\label{A1}
If $\deg L$ is odd, then $A_{3g-4}(\oSU)\simeq\bZ$.
\end{cor}

\pf
Let $\deg L=2g-1$.
The isomorphism $\phi_{L,3}|_U\colon U\fun V$ of Proposition \ref{UVproper}
induces an isomorphism $(\phi_{L,3}|_U)_*\colon A_{3g-4}(U)\fun A_{3g-4}(V)$,
and $A_{3g-4}(V)\simeq A_{3g-4}(\oSU)$ because the complement of $V$ has
codimension $\geq2$ in $\oSU$.
Recall that $U$ is an open subset of $\bP_{L,3}$ whose complement has
codimension one, and therefore there exists an exact sequence (see [Ful84,
1.8])
$$ A_{3g-4}(\bP_{L,3} \setminus U) \lfun A_{3g-4}(\bP_{L,3}) \lfun A_{3g-4}(U)
\lfun 0, $$
where
$$ A_{3g-4}(\bP_{L,3}) \simeq \bZ H \oplus_{p \in J} \bZ \tE_1|_p
\oplus_{p \in J} \bZ \tE_{2,p} \oplus \bZ E_3, $$
with $H$ the pull-back of a hyperplane class from $\bP_L$.
It follows from our description of $U$ in the proof of Proposition
\ref{UVproper} that
$$ A_{3g-4}(\bP_{L,3} \setminus U) \simeq \bZ \tG_L \oplus_{p \in J}
\bZ \tE_1|_p \oplus_{p \in J} \bZ \tE_{2,p}, $$
and therefore
$$ A_{3g-4}(U) \simeq \frac{\bZ H \oplus \bZ E_3}{\bZ \tG_L}. $$
Since, as in the smooth case, $\tG_L\sim gH-(g-1)E_3$, we obtain
$$ A_{3g-4}(\oSU) \simeq \frac{\bZ H \oplus \bZ E_3}{\bZ (g H - (g - 1)E_3)}
\simeq \bZ. $$
\qed

We now study the complement of $\SU_C(2,L)$ in $\oSU$.

\begin{pro}\label{codimSUoSU}
For every irreducible nodal curve $C$ of genus $\geq2$, 
$$ \emph{codim}(\oSU \setminus \SU_C(2,L), \oSU) \geq 3. $$
\end{pro}

\pf
If suffices to prove this in the case when $\deg L=d>>0$.
The generic element of $\oSU\setminus\SU_C(2,L)$ is a torsion-free
non-locally-free sheaf such that every section vanishes at a node.
Therefore, it is a non-locally-free extension of the form
$$ 0 \lfun (\pi_p)_* \cO_{C_p} \lfun E \lfun L \otimes ((\pi_p)_* \cO_{C_p})^*
\lfun 0 $$
for some node $p$ in $C$, i.e., the push-forward of an extension
$$ 0 \lfun \cO_{C_p} \lfun \cE \lfun \pi_p^*L (- p_1 - p_2) \lfun 0 $$
via the partial normalization $\pi_p\colon C_p\fun C$.
A generic such extension $\cE$ is in $\SU_{C_p}(2,\pi_p^*L(-p_1-p_2))$, and
there exists a morphism
$$ \ov{\SU_{C_p}(2, \pi_p^*L (- p_1 - p_2))} \overset{(\pi_p)_*}{\lfun} \oSU$$
defined by $\cE\mapsto(\pi_p)_*\cE$.
It is a morphism because if $F\su(\pi_p)_*\cE$ is a rank-$1$ torsion-free
subsheaf of $(\pi_p)_*\cE$, then $\pi^*F/\Tors\su\cE$.
Since $\cE$ is semi-stable, the degree of $\pi^*F/\Tors$ is
$\leq\deg\cE/2=(d-2)/2$.
Therefore $\deg F\leq(d-2)/2+1=d/2$, and $(\pi_p)_*\cE$ is semi-stable.

The dimension of $\oSU\setminus\SU_C(2,L)$ is therefore less than or equal to
the dimension of $\SU_{C_p}(2,\pi_p^*L(-p_1-p_2))$, which is $3g-6$ unless
$g=2$ and $\deg L$ is even when it is $1=3g-5$.
But if $g=2$ and $\deg L$ is even then $\oSU=\SU_C(2,L)\simeq\bP^3$ (see
[Bho98]).
Therefore,
$$ \codim(\oSU \setminus \SU_C(2,L), \oSU) \geq (3g - 3) - (3g - 6) = 3. $$
\qed

\begin{rem}
The proof of the proposition shows that, except for the case $g=2$ and
$\deg L$ even, $\oSU\setminus\SU_C(2,L)$ is the union of $|J|$ varieties of
dimension $3g-6$, the images of the morphisms $(\pi_p)_*$ described above, one
for each node.
\end{rem}

An immediate consequence of Corollary \ref{A1} and Proposition
\ref{codimSUoSU} is the following result, which was already proved by Bhosle
(see [Bho99] and [Bho04]).

\begin{cor}
If $C$ is an irreducible nodal curve of genus $\geq2$ and $\deg L$ is odd,
then
$$ A_{3g-4}(\SU_C(2,L)) \simeq \bZ. $$
\end{cor}

Our last application is the following corollary.

\begin{cor}
If $C$ is an irreducible nodal curve of genus $\geq2$, then
$A_{3g-4}(\oSU)\simeq\bZ$.
\end{cor}

\pf
Bhosle proved in [Bho99] and [Bho04] that $A_{3g-4}(\SU_C(2,L))\simeq\bZ$.
The result follows from the exact sequence (see [Ful84, 1.8])
$$ A_{3g-4}(\oSU \setminus \SU_C(2,L)) \lfun A_{3g-4}(\oSU) \lfun
A_{3g-4}(\SU_C(2,L)) \lfun 0 $$
and Proposition \ref{codimSUoSU}.
\qed

\end{document}